\documentclass[journal,twoside,web]{ieeecolor}

\usepackage{generic}
\usepackage{cite}
\usepackage{amsmath,amssymb,amsfonts}
\usepackage{algorithmic}
\usepackage{graphicx}
\usepackage{textcomp}

\usepackage[colorlinks=true, allcolors=cyan]{hyperref}

\usepackage{upgreek}
\usepackage{nicefrac}
\usepackage[para,online,flushleft]{threeparttable}

\newcommand{\Sabs}[1]{\left\lceil #1 \right\rfloor}

\newcommand{\abs}[1]{\left| #1 \right|}
\newcommand{\norm}[1]{\left|\left| #1 \right|\right|}

\newcommand{\conjunto}[1]{\left\{ #1 \right\} }

\newcommand{\Subsin}[1]{\left. #1 \right|}

\newtheorem{Lema}{Lemma}
\newtheorem{theorem}{Theorem}
\newtheorem{remark}{Remark}
\newtheorem{definition}{Definition}

\def\BibTeX{{\rm B\kern-.05em{\sc i\kern-.025em b}\kern-.08em
    T\kern-.1667em\lower.7ex\hbox{E}\kern-.125emX}}
\begin{document}

\title{Earthquake Control: An Emerging Application for Robust Control. Theory and Experimental Tests}

\author{Diego~Guti\'errez-Oribio, Georgios Tzortzopoulos, Ioannis Stefanou and Franck Plestan
\thanks{Diego~Guti\'errez-Oribio, Georgios Tzortzopoulos and Ioannis Stefanou are in the \'Ecole Centrale de Nantes - GeM, 1 Rue de la No\"e, 44321 Nantes, France (e-mail: diego.gutierrez-oribio@ec-nantes.fr, georgios.tzortzopoulos@ec-nantes.fr, ioannis.stefanou@ec-nantes.fr). }
\thanks{Franck Plestan is in the \'Ecole Centrale de Nantes - LS2N, UMR CNRS 6004, 1 Rue de la No\"e, 44321 Nantes, France (e-mail: franck.plestan@ec-nantes.fr).}}

\maketitle

\begin{abstract}
This paper addresses the possibility of using robust control theory for preventing earthquakes through fluid injections in the earth's crust. The designed robust controllers drive aseismically a fault system to a new equilibrium point of lower energy by tracking a slow reference signal. The control design is based on a reduced-order nonlinear model able to reproduce earthquake-like instabilities. Uncertainties related to the frictional and mechanical properties of the underlying physical process and external perturbations are considered. Two types of controllers are derived. The first one is based on sliding-mode theory and leads to local finite-time convergence of the tracking error and rejection of Lipschitz w.r.t. time perturbations. The second controller is based on LQR control and presents global exponential stability of the tracking error and rejection of Lipschitz w.r.t. states perturbations. Both controllers generate a continuous control signal, attenuating the chattering effect in the case of the sliding-mode algorithms. The developed controllers are tested extensively and compared on the basis of numerical simulations and experiments in the laboratory. The present work opens new perspectives for the application of robust nonlinear control theory to complex geosystems, earthquakes and sustainable energy production.
\end{abstract}


\section{Introduction}
\label{sec:Introduction}


\IEEEPARstart{E}{}arthquakes are dynamic instabilities that occur in the earth's crust. $65\%$ of the most catastrophic earthquakes happen at depth up to $12$ km and are of natural causes \cite{b:Scholz-2002}. However, earthquakes also occur due to anthropogenic causes. It is nowadays established that injecting fluids in the earth's crust can reactivate seismic faults, leading to important earthquake events (see \cite{b:Rubinstein-Mahani-2015}, \cite{b:Keranen-Savage-Abers-Cochran-2013} and \cite{b:Zastrow-2019}, to name a few examples).

In this paper, fluid injections are seen from another perspective. Instead of considering them as an earthquake triggering mechanism, they are seen as an input to a dynamical system, which can stabilize it and achieve tracking over a reference signal, if it is adequately designed. This dynamical system is the physical process leading to earthquake instabilities and it is characterized by important nonlinearities due to friction. Moreover, it can present many uncertainties and unmodelled dynamics, that a successful control scheme has to compensate. Finally, unlike many existing applications of control theory that target decreasing the response time of the system, here the aim is the opposite, \textit{i.e.}, to slow down the system dynamics. The above mentioned characteristics of this system result in a challenging problem for control theory, with important applications in energy production (\textit{e.g.}, oil, gas, deep geothermal energy, and CO$_2$ sequestration) and earthquake prevention.



Existing strategies for earthquake control are very limited and ad-hoc. One can refer, for instance, to the field experiments in Rangely, Colorado, US, \cite{b:Raleigh1-Healy1-Bredehoeft}, where seismicity was reduced by turning off the pore pressure. In Dale, New York, US, where earthquakes of magnitude $1$ to $1.4$ were arrested by dropping the top hole pressure below $5$ MPa \cite{b:Fletcher-Sykes}. More recent field experiments involve the well monitored tests by \cite{b:Guglielmi-Cappa-Avouac-Henry-Elsworth,b:Cappa-Scuderi-Collettini-Guglielmi-Avouac}. However, as mentioned above, all these experiments were based on trial and error and they were not based on control theory.

Recently, an LQR control was designed to stabilize and perform tracking of an earthquake modelled by a MIMO system \cite{b:Stefanou2019}, whereas a double-scale asymptotic approach was employed to design a transfer function-based control in \cite{b:Stefanou2020}. These first applications of control theory to this problem have shown that earthquakes could be controlled, at least from a mathematical point of view.

The objective of this paper is twofold. First, to evaluate the performance of different controllers and, second, to test them in the laboratory with a specially designed apparatus \cite{b:Stefanou-Tzortzopoulos-Braun-Patent}. The design of the controllers is based on a reduced model for earthquakes. This reduced-order model establishes an average behaviour of a single earthquake fault (see \cite{b:Scholz-2002}, among others). It consists of a single mass that can slide on a rough surface under friction. The frictional interface is usually a complex structure (see \cite{b:Ben-Zion-Sammis-2003,b:Brodie-Fettes-Harte-Schmid-2007}), where various physico-chemical phenomena take place during seismic slip (see \cite{b:Scholz-2002,b:Rattez-Stefanou-Sulem-2018,b:Rattez-Stefanou-Sulem-Veveakis-Poulet-2018,b:Reches-Lockner-2010}). As a result, the friction coefficient and its weakening, not only depend on the slip and the slip-rate, but also on the evolution of the microstructural network, the grain size, the presence and pressure of interstitial fluids, the temperature, time, the reactivation of chemical reactions and other multiphysics couplings (see \cite{b:Brantut-Sulem-2012,b:Veveakis-Alevizos-Vardoulakis-2010,b:Veveakis-Stefanou-Sulem-2013,b:Collins-craft2020}).
 
All these complex phenomena induce the presence of uncertainties and/or perturbations to the plant. As a consequence, they need to be compensated by a robust controller able to obtain a slow-aseismic response. A classic robust approach versus constant perturbations, is the integral action (see \cite[Chapter 12]{b:Khalil2002}). Among many robust control approaches, one can also cite the sliding-mode theory \cite{b:utkin92,b:SMC_Fridman}. This type of control is known for being insensitive to bounded and matched perturbations, leading to finite-time convergence. Sliding-mode controllers are also known for the simplicity of their gains tuning. The problem with these controllers though is the use of a discontinuous function, \textit{sign}, that may lead to the so-called \textit{chattering} effect, possibly damaging the actuators. Recently, in order to address the above mentioned drawback, the Continuous Higher-Order Sliding-Modes Algorithms (CHOSMA, see \cite{b:Chalanga-Kamal-Bandyopadhyay,b:Torres-Sanchez-Fridman-Moreno,b:Moreno_2016}) have been developed to keep the interesting features of the sliding-mode algorithms (compensation of Lipschitz w.r.t. time perturbations in finite-time), while using a continuous control signal.

In this paper, two control strategies are proposed to achieve a slow-seismic response: the first one is based on CHOSMA and the other based on a Linear Quadratic Regulator extended with integral action (e-LQR based on the original LQR control in \cite{b:Kalman-1960}). The former is robust against Lipschitz w.r.t. time perturbations, while the latter is robust against Lipschitz w.r.t. state perturbations, presenting both continuous control signal. The designed controllers stabilize and reduce the natural response time of the system, making the energy dissipation to be many orders of magnitude slower compared to a real earthquake event. This is an uncommon paradigm in control theory, where usually the objective is to drive the states of the system to the origin as fast as possible. Finally, to test the feasibility of the presented algorithms, simulations and experimental confirmation using a real laboratory test benchmark, able to reproduce earthquake-like instabilities, are presented. This allows testing the controllers and comparing their performance.

The outline of this work is as follows. The description of the reduced-order model  for reproducing earthquake-like instabilities, its instability condition, and the control objectives are given in Section \ref{sec:Problem}. The experimental setup is shown in the same section, while the design of the \textit{robust} control strategies is detailed in Section \ref{sec:Control}. Simulations and experimental results are shown in Section \ref{sec:SimExp}, and concluding remarks are made in Section \ref{sec:Conclusions}. 

\subsection{Preliminaries}

Throughout the text, the term ''Lipschitz w.r.t. the time/state function" is used to call a function that fulfils a Lipschitz condition with respect to the time/state. Let us consider the time-varying differential equation
    \begin{equation*}
    \dot{x}=f(t,x(t)),\; t\ge t_{0},
    \end{equation*}
where $x(t)\in{\mathbb{R}}^{n}$ is the state vector; $f:{\mathbb{R}}_{\geq 0}\times{\mathbb{R}}^{n}\to{\mathbb{R}}^{n}$ is a function that can be discontinuous, measurable with respect to $t$, and $f(t,0)=0$. The initial condition $x(t_{0})\in{\mathbb{R}}^{n}$ at time instant $t_{0}\in\mathbb{R}$ is denoted as $x_0$. The solution of the system is understood in the Filippov's sense (see \cite{b:filippov}). Let $\Omega$ be an open subset of ${\mathbb{R}}^{n}$, such that $0\in\Omega$. 

\begin{definition}
\label{def:stability}
\cite{b:Khalil2002,b:Bhat-Bernstein}.
The origin, $x = 0$, of the latter system is said to be:\\
$\bullet$ Locally Stable (LS) if for any $\epsilon >0$ there is $\delta=\delta(\epsilon,t_0)>0$ such that if $\norm{{ x }_{ 0 }} \le \delta$ then $\norm{x(t)} \le \epsilon$ for any ${x}_{0} \in \Omega$ and for all $t\ge { t }_{ 0 }$.\\
$\bullet$ Locally Asymptotically Stable (LAS) if it is LS and there is $c=c(t_0) > 0$ such that if $\norm{{ x }_{ 0 }} \le c$ then $x(t) \rightarrow 0$ for any ${x}_{0} \in \Omega$ and for all $t\ge { t }_{ 0 }$.\\
$\bullet$ Locally Exponentially Stable (LES) if it is LAS and there are $k>0$, $\lambda>0$ such that $\norm{x(t)} \leq k\norm{x_0}e^{-\lambda(t-t_0)}$ for any ${x}_{0} \in \Omega$ and for all $t\ge { t }_{ 0 }$.\\
$\bullet$ Locally finite-time Stable (LFTS) if it is LAS and $x(t)=0$ for all $t\ge T(t_0,{ x }_{ 0 } )$, where $T:\mathbb{R}_{\geq 0} \times \mathbb{R}^n \to \mathbb{R}_{\geq 0}$ is called the settling-time function.\\
If $\Omega={\mathbb{R}}^{n}$, then $x = 0$ is said to be Globally Stable (GS), Globally Asymptotically Stable (GAS), Globally Exponentially Stable (GES), Globally finite-time Stable (GFTS), respectively.
\end{definition}

The definition of weighted homogeneity is introduced to be used in the sequel.
\begin{definition} \cite{b:Baccioti-Rosier_2005}, \cite{b:Bernuau-Efimov-Perruquetti-Polyakov}. 
Consider the vector $x \in \mathbb{R}^n$. Its dilation operator is defined as $\Delta_\epsilon^{r}x: = (\epsilon^{r_1} x_1, ..., \epsilon^{r_n} x_n)$, $\forall \epsilon > 0$, where $r_i > 0$ are the weights of the coordinates and $r = (r_1, ..., r_n)$ is the vector of weights. A function $V: \mathbb{R}^n \rightarrow \mathbb{R}$ (or a vector field $f: \mathbb{R}^n \rightarrow \mathbb{R}^n$, or vector-set $F(x) \subset \mathbb{R}^n$) is called $r$-homogeneous of degree $m \in \mathbb{R}$ if the identity $V(\Delta_\epsilon^{r}) = \epsilon^{m}V(x)$ holds (or $f(\Delta_\epsilon^{r} x) = \epsilon^{m} \Delta_\epsilon^{r} f(x)$, or $F(\Delta_\epsilon^{r} x) = \epsilon^{m} \Delta_\epsilon^{r} F(x)$). 
\end{definition}

The following result is well-known for continuous homogeneous functions (see \cite{b:Andrieu-Praly-Astolfi_2008, b:Moreno_2016}), and can be extended to semi-continuous functions \cite{b:Cruz-Zavala-Moreno_2017}:
\begin{Lema}
	Let $\eta: \mathbb{R}^n \rightarrow \mathbb{R}$ and $\gamma: \mathbb{R}^n \rightarrow \mathbb{R}$ be two $r$-homogeneous and upper semi-continuous single-valued functions, with the same weights $r = (r_1, ..., r_n)$ and homogeneity degree $m>0$. Suppose that $\gamma(x) \leq 0$ in $\mathbb{R}^n$. %
	 If
\begin{equation*}
  \{x \in \mathbb{R}^n \setminus \{0\}: \gamma(x) = 0 \} \subseteq \{x \in \mathbb{R}^n \setminus \{0\}: \eta(x) < 0 \},
\end{equation*}
then there exists a real number  $\lambda^*$ and a constant $c > 0$ so that, for all $\lambda \geq  \lambda^*$ and for all $x \in \mathbb{R}^n \setminus\{0\}$ the following inequality is satisfied $\eta(x) + \lambda \gamma(x) \leq -c \norm{x}_{r, p}^m$. 
	\label{lem:Negative I Semi}
\end{Lema}

Define the function $\lceil\cdot\rfloor^{\gamma}:=|\cdot|^{\gamma}\mathrm{sign}(\cdot)$, for any $\gamma\in \mathbb{R}_{\geq 0}$ with $\mathrm{sign}(x)=\left\{\begin{array}{cc}
1 & x>0 \\ 
\left[-1,1\right] & x=0 \\ 
-1 & x<0
\end{array}   \right.$.

The following Lemma is simple (just monotonicity) but useful
\begin{Lema}
	\cite{b:Mercado-Uribe-Moreno} 
	Consider the real variables $x$, $y$, it is always true that $\textup{sign } \big(\Sabs{x + y}^\beta - \Sabs{y}^\beta \big) = \textup{sign }(x), \quad \beta > 0$.
	\label{lem:Holder S}
\end{Lema}

\section{Problem Statement}
\label{sec:Problem}

Two scenarios for earthquake modelling used in this paper are presented in this Section. The first one is a reduced model of a real earthquake on which the control designs are performed and numerical simulations are made, whereas the second one, is a novel experimental setup designed to reproduce and then control earthquake-like instabilities in the laboratory.  

\subsection{Reduced Model for Earthquakes}

\begin{figure}[ht!]
  \centering 
  \includegraphics[width=8.7cm,height=3.5cm]{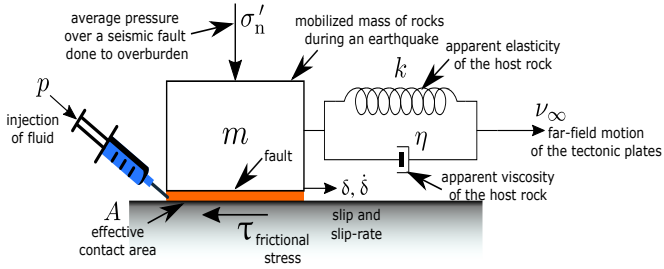}
  \caption{Reduced mechanical model for reproducing earthquake-like instabilities.}
  \label{fig:spring-slider}
\end{figure}

The dynamics of earthquakes can be represented, in average/energetical sense, with the spring-slider analogue system (see \cite{b:Stefanou2019,b:Stefanou2020,b:Scholz-2002,b:Kanamori-Brodsky-2004,b:Tzortzopoulos-Braun-Stefanou-2021}) depicted in Fig. \ref{fig:spring-slider}. 

This mechanical system consists of a mass, $m$, which slides on a frictional interface (equivalent to a seismic fault). The mass is connected to a Kelvin-Voigt configuration composed by a spring with stiffness $k$ (equivalent to the apparent elasticity of the host rock) and a dashpot with damping coefficient $\eta$ (equivalent to the apparent viscosity of the host rock). At the other extremity of the Kelvin-Voigt configuration a constant velocity, $\nu_\infty$, is applied (equivalent to the far-field motion of the tectonic plates). It is assumed Coulomb friction with a friction coefficient $\mu(\delta,\dot{\delta})$ that depends on the slip $\delta$ (block's displacement) and the slip-rate $\dot{\delta}$ (block's velocity). According to Terzaghi's principle, the frictional stress $\uptau$ takes the following form: $\uptau=\mu(\delta,\dot{\delta})\hspace{2pt}(\sigma^{\, \prime}_{\text{n}}-p)$, where $\sigma^{\, \prime}_{\text{n}}$ is the constant/reference average effective stress (\textit{e.g.}, the overburden due to the weight of the rocks and the interstitial fluid pressure) and $p$ the fluid pressure developed due to injecting fluid. $p$ is the input to the system for which the controllers will be designed and tested.


According to \cite{b:Scholz-2002} and \cite{b:Kanamori-Brodsky-2004}, approximately a rock mass of volume $L_{\text{ac}}^3$ is mobilized during an earthquake event, where $L_{\text{ac}}$ is equal to the length of the seismic fault. Therefore, the mobilized mass during an earthquake event is $m \approx \rho L_{\text{ac}}^3$, where $\rho$ is the density of the surrounding to a seismic fault rocks. The fault length can be calculated as $L_{\text{ac}}=\nicefrac{G}{\bar{k}}$, where $G$ is the shear-modulus of the host rock and $\bar{k} = \nicefrac{k}{L_{\text{ac}}^2}$, its apparent normalized elastic stiffness. The damping coefficient $\eta$ is given by $\eta=2\zeta m \omega_n$, where $\zeta$ is the damping ratio and $\omega_n=\sqrt{\nicefrac{k}{m}}$, the natural frequency of the reduced system.

Applying the force balance equation, system in Fig. \ref{fig:spring-slider} can be represented by the following mathematical model
\begin{equation}
m \ddot{\delta} = -\mu(\delta,\dot{\delta})A(\sigma^{\, \prime}_{\text{n}}-p) + k(\delta_\infty-\delta) + \eta(\nu_\infty-\dot{\delta})+\varphi_e(\delta,\dot{\delta},t),
\label{eq:springmodel}
\end{equation}
where $A \approx L_{\text{ac}}^2$ is the effective contact area (fault rupture area), $\delta_\infty=\nu_\infty t$ the displacement at the extremity of the Kelvin-Voigt configuration, $\ddot{\delta}$ the acceleration of the mobilized block, and $\varphi_e(\delta,\dot{\delta},t)$ is a perturbation affecting the system, \textit{e.g.}, an external perturbation or unmodelled dynamics due to the complex frictional phenomena. 

In this paper, the friction coefficient $\mu(\delta,\dot{\delta})$ is assumed to fulfil
\begin{equation}
  0<\mu_{res} \leq \mu(\delta,\dot{\delta}) \leq \mu_{max},
  \label{eq:mug}
\end{equation}
where the constants $\mu_{res}$ and $\mu_{max}$ are given. Such assumption is fulfilled by friction laws used in fault mechanics, like the slip-weakening friction law \cite{b:Kanamori-Brodsky-2004}, the slip-rate weakening law \cite{b:Huang1992a}, and the rate-and-state friction law \cite{b:Dieterich1981,b:Ruina-1983}. Notice that the exact frictional rheology is not known in reality, which needs the design of robust controllers, as will be done in Section \ref{sec:Control}.

Additionally, according to \cite{b:Kanamori-Brodsky-2004}, the seismic magnitude $M_w$ is defined as
\begin{equation}
  M_w = \frac{2}{3} \log_{10}M_0-6.07, \quad M_0 = L_{ac}^3 \Delta \uptau,
  \label{eq:magnitude}
\end{equation}
where $M_0$ is the seismic moment measured in [Nm] and $\Delta \uptau=(\mu_{max}-\mu_{res})\sigma^{\, \prime}_{\text{n}}$.

\subsection{Shifted System and Instability Condition} 

Defining the state variables $z_1=\delta$ and $z_2=\dot{\delta}$, the state representation of system \eqref{eq:springmodel} is 
\begin{equation*}
\begin{split}
  \dot{z}_1&=z_2,\\
  \dot{z}_2&=-\mu(z_1,z_2)\hat{N}(\sigma_n^\prime-p) + \hat{k}(\delta_\infty-z_1) + \hat{\eta}(v_\infty-z_2) \\ &\quad+\hat{\varphi}_e(z_1,z_2,t),
\end{split}  
\end{equation*}
where $\hat{N}=\nicefrac{A}{m}$, $\hat{k}=\nicefrac{k}{m}$, $\hat{\eta}=\nicefrac{\eta}{m}$ and $\hat{\varphi}_e(z_1,z_2,t)=\nicefrac{\varphi_e(z_1,z_2,t)}{m}$.

The set of equilibrium points $(z_1^*,z_2^*)$ of the above system in open loop and without perturbation, \textit{i.e.} $\hat{\varphi}_e(z_1,z_2,t)=0$, is described by
\begin{equation*}
\begin{split}
  z_1^* = -\mu(z_1^*,z_2^*)\frac{\hat{N}}{\hat{k}}\sigma_n^\prime + \delta_\infty + \frac{\hat{\eta}}{\hat{k}}v_\infty, \quad
  z_2^* = 0.
\end{split}
\end{equation*}
Note that the equilibrium $(z_1^*,z_2^*)$ depends on the friction coefficient $\mu(z_1^*,z_2^*)$. In this paper, the controller design will be made when the system reaches the above equilibria, which can be unstable. The maximum value of the friction coefficient, $\mu(0,0) = \mu_{max}$, is considered for being on the verge of slip and the system is shifted as follows. Setting $z_1^*=0$ and the new state variables $x_1 = z_1-z_1^*$ and $x_2 = z_2-z_2^*$, the shifted system reads as
\begin{equation}
\begin{split}
  \dot{x}_1&=x_2,\\
  \dot{x}_2&=-[\mu(x_1,x_2)-\mu^*]\hat{N}\sigma_n^\prime + \mu(x_1,x_2)\hat{N} p -\hat{k}x_1 - \hat{\eta}x_2 \\ &\quad + \hat{\varphi}_e(x_1,x_2,t),
\end{split}  
\label{eq:shift}
\end{equation}
where $\mu^* = \mu(0,0) = \mu_{max}$. Note that if $\hat{\varphi}_e(x_1,x_2,t)=0$, system \eqref{eq:shift} has an equilibrium point located at the origin $x_1^*=x_2^*=0$ in open loop. 

In order to analyse the stability of the origin of system \eqref{eq:shift} without the perturbation term $\hat{\varphi}_e(x_1,x_2,t)$, consider its Jacobian matrix $J(x_1,x_2)$ evaluated at the origin as
{\small
\begin{equation*}
\begin{split}
  &J(0,0) = \\
  &\left[\begin{array}{cc}
  0 & 1 \\ 
  -\hat{k}-\hat{N}\sigma_n^\prime \Subsin{\frac{\partial \mu}{\partial x_1}}_{(x_1,x_2)=(0,0)} & -\hat{\eta}-\Subsin{\frac{\partial \mu}{\partial x_2}}_{(x_1,x_2)=(0,0)}
  \end{array} \right],
\end{split}
\end{equation*}}
where the conditions to have an unstable origin are
\begin{equation}
\begin{split}
k < -A\sigma_n^\prime \Subsin{\frac{\partial \mu}{\partial x_1}}_{(x_1,x_2)=(0,0)}, \quad
\eta < -m \Subsin{\frac{\partial \mu}{\partial x_2}}_{(x_1,x_2)=(0,0)}.
\end{split}
\label{eq:inst_cond}
\end{equation}

The above inequalities are in accordance with the nominal studies of \cite{b:Stefanou2019,b:Scholz-2002,b:Dieterich1979} and they show that dynamic instability will take place when the elastic unloading of the springs or the apparent viscosity of the host rock cannot be counterbalanced by friction. 

\subsection{Experimental setup}

\begin{figure*}[ht!]
  \centering 
  \includegraphics[width=13cm,height=7cm]{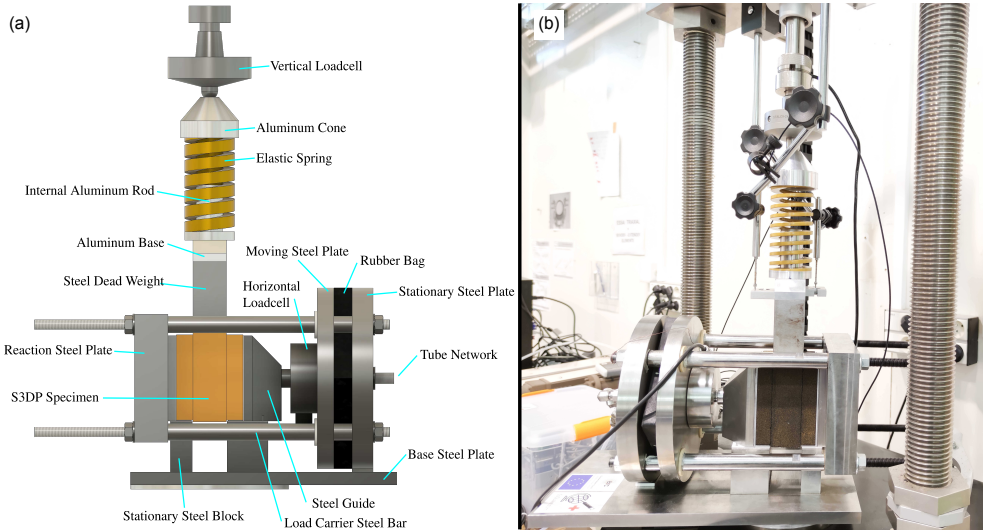}
  \caption{Experimental apparatus for reproducing and controlling earthquake-like instabilities: (a)Schematic figure, (b)Real configuration}
  \label{fig:spring-slider-lab}
\end{figure*}

A novel experimental apparatus for reproducing and controlling earthquake-like instabilities in the laboratory, has been designed (see \cite{b:Tzortzopoulos-2021,b:https://doi.org/10.1029/2021JB023410,b:Stefanou-Tzortzopoulos-Braun-Patent}). This experimental setup is depicted in Fig. \ref{fig:spring-slider-lab}.

Two loading systems are used in order to apply shear (displacement controlled by the vertical loading system) and normal (pressure controlled by the horizontal loading system) stresses to the sheared interfaces. The horizontal loading system consists of an inflatable rubber bag whose pressure can be controlled, in real-time, through a fast response electro-pneumatic pressure regulator. This system can simulate fluid injection/extraction into/from the fault interface by properly adjusting the (air) pressure in the rubber bag to desired levels resulting in variations of the effective normal stress in the sheared interfaces. The vertical loading system consists of a press compressing slowly a linear-elastic spring, which simulates the stored energy of the earthquake. 

The specimen consists of three 3D printed samples of sand particles (see \cite{b:Braun-Tzortzopoulos-Stefanou-2021} for more details about this surrogate material and its characterization) and it is located below the spring. The middle block of the specimen simulates the mobilized mass of the rocks and its frictional interfaces with the adjacent blocks the seismic fault. No specific elements are used for damping, which is provided naturally by the various components of the experimental setup and the specimen itself. The friction coefficient between the blocks can be measured in the experimental setup (see \cite{b:Tzortzopoulos-2021,b:Braun-Tzortzopoulos-Stefanou-2021}). However, these tests are not presented here because the objective is to design robust controllers that are agnostic to the frictional rheology. The frictional rheology in a real fault might depend on the slip, slip-rate and time, but also in thermal, chemical and other processes. Moreover, due to heterogeneity, the frictional properties of faults involve uncertainties and contrary to a laboratory experiment, they cannot be inferred with high accuracy. Therefore, successful controllers must be able to compensate all these sources of uncertainties and to be able to obtain slow aseismic response with minimum information about the system.

For measuring the applied vertical load and the horizontal forces, load cells are used. For measuring the average slip of the middle block, two vertically placed Linear Variable Differential Transformers (LVDTs) are used. Their readings are averaged to eliminate parasitic measurements related to possible rotation of the middle block. Finally, the pump for controlling the pressure in the rubber bag simulates the fluid injection to the fault and can supply pressures up to 1 [MPa]. Note that in this setup, increasing the pressure in the rubber bag corresponds to decreasing the pressure, $p$, in the fault (see also eq. \eqref{eq:springmodel}). However in the following, the distinction between pressure in the bag and pressure in the fault is not made and the results are presented in function of the pressure at the fault (\textit{i.e.}, $p=-p_{bag}$). All these sensors are connected to a data acquisition device processed using the LabVIEW software, which processes the data with a sampling rate of $1$ [ms]. LabVIEW also allows the implementation of our robust controllers for controlling earthquake-like instabilities. See \cite[Chapter 5]{b:Tzortzopoulos-2021} for more details about the experimental setup.

\begin{remark}
The electro-pneumatic pressure regulator has a response rate of $1$ [ms]. It is controlled by a PID algorithm, which is faster than the response rate of the pump. The operating system of the computer in which the controllers were implemented is of $1$ to $2$ [ms]. The above characteristics times are much lower than the characteristic time of the system instability, which is $50$ [ms].Therefore, the actuator dynamics (the pump) was not taken into account in the system \eqref{eq:springmodel}.
\end{remark}

\subsection{Control Objective}

\begin{figure}[ht!]
  \centering 
  \includegraphics[width=8cm,height=4cm]{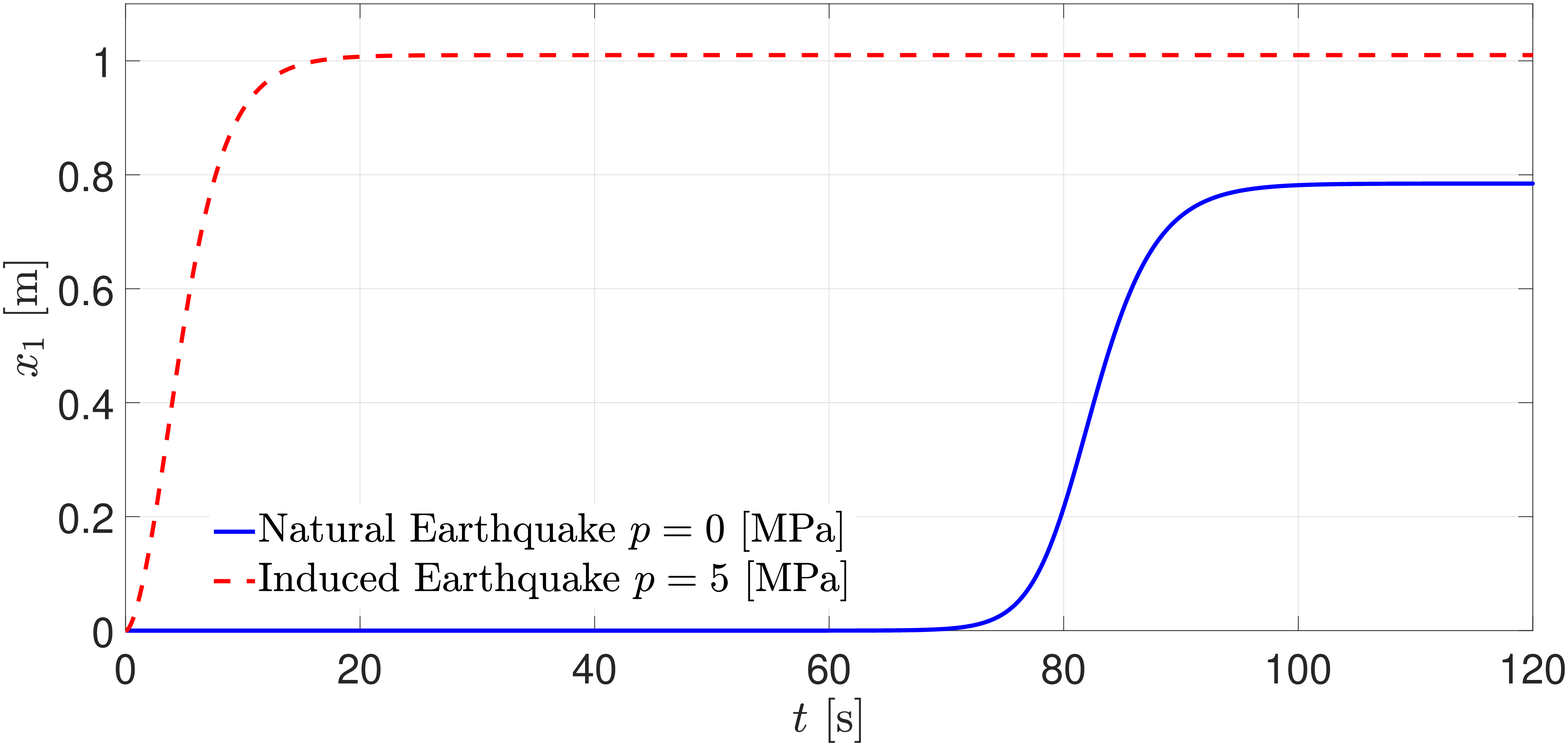}
  \includegraphics[width=8cm,height=4cm]{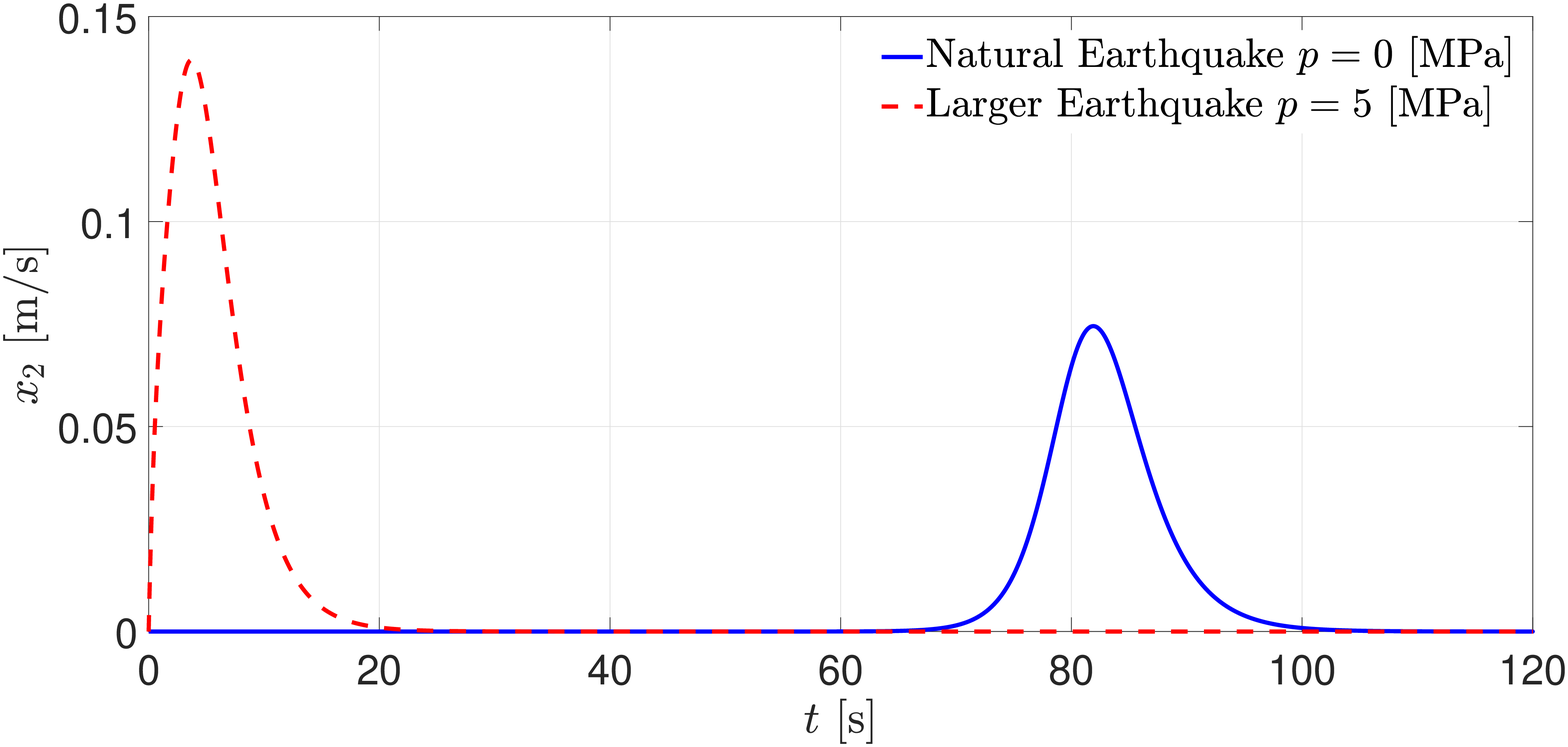}
  \caption{Slip and slip-rate in terms of time for two different scenarios. Blue curve: Natural earthquake. Red Dashed Curve: Induced earthquake.}
  \label{fig:sim}
\end{figure}

As shown in system \eqref{eq:shift} and Fig. \ref{fig:spring-slider}, the fluid pressure $p$ is the only input acting on the dynamics of the mechanical system. In a real-scale scenario, fluid injections in the earth's crust change the fluid pore pressure over seismic faults \cite{b:Cappa-Scuderi-Collettini-Guglielmi-Avouac}. As shown in \cite{b:Tzortzopoulos-Braun-Stefanou-2021} among others, this can destabilize the fault system and induce/trigger larger earthquakes. In order to illustrate this phenomenon, two numerical simulations of the dynamical system \eqref{eq:shift} are shown in Fig. \ref{fig:sim}. The first one has no input pressure (Natural earthquake using $p=0$ [MPa]) and the second one has a constant pressure (Induced earthquake using $p=5$ [MPa]). Both simulations were made starting from the origin and with an external perturbation $\hat{\varphi}_e=3.2\times 10^{-15}$ [m/s$^2$] applied at $t=0$ [s] just to move the states out of the equilibrium point. Note how an earthquake event ($\delta$ evolving with high rate) has been triggered in both cases, but a larger and faster one was obtained on the second case where $p>0$. The system parameters used for these simulations are shown in Table \ref{tab:para} (Real-fault column).

In the previous paragraph, the fact that open-loop fluid injections in the earth's crust might risk to stimulate large seismic events is highlighted. To prevent this, one could adequately adjust the fluid pressure (input $p$) by employing control techniques in order to stabilize the system \eqref{eq:shift} origin and/or track a reference input signal, releasing the stored elastic energy smoothly ($\delta$ evolving slowly) and not abruptly as shown above. 

From a mathematical point of view, a reduced-order model for earthquakes like the system \eqref{eq:springmodel} captures well the low-frequency instabilities, which are of most importance, because they carry most of the seismic energy. For the experimental confirmation, the designed apparatus in Fig. \ref{fig:spring-slider-lab} reproduces these low-frequency instabilities and successfully designed controllers must be able to dissipate them slowly. See \cite{b:https://doi.org/10.1029/2021JB023410,b:Gutierrez-Stefanou-Plestan-2022}, for theoretical works for the control of more detailed models of the physical process of rupture of a single mature seismic fault. 

Therefore, the objective in the Section \ref{sec:Control} is to design a control law $p$ driving $\delta$ and $\dot{\delta}$ in system \eqref{eq:shift} to follow some desired predefined references of slow slip rate, $r(t),\dot{r}(t)$, resulting in a slow-aseismic response, despite the presence of $\hat{\varphi}_e$ and parametric uncertainties. Furthermore, this control design will be tested in Section \ref{sec:SimExp} on simulations and on the experimental setup depicted in Fig. \ref{fig:spring-slider-lab}.

The designed controllers will reduce the natural response time of the system slowing its energy dissipation and eliminating bursts of kinetic energy (earthquake phenomenon). This is an uncommon paradigm in control theory where usually the objective is to drive the states of the system to the origin as fast as possible

\section{Control Design}
\label{sec:Control}

Two control approaches are presented in this section: the first one is based on continuous sliding-mode control theory \cite{b:Torres-Sanchez-Fridman-Moreno,b:Moreno_2016} whereas the second one is based on LQR control \cite{b:Kalman-1960}. Both these controllers must be able to force the states of system \eqref{eq:shift} to track a predefined reference, using a continuous control signal, in spite of uncertainties/perturbations of the system \eqref{eq:shift}.

It is worth emphasizing that due to the nature of the physical process of the earthquake phenomenon, the controllers have to be robust to uncertainties regarding the material properties and plant dynamics, as well as to disturbances of the input signals and measurements. This motivates the choice of sliding-modes control, which will be compared with a robust LQR control in order to slow down the system and make it follow a desired reference trajectory.

The desired reference for the output $y=x_1$ is a smooth function reading as
\begin{equation}
    r(t) = d_{max}s^3(10-15s+6s^2),
  \label{eq:ref}
\end{equation}
where $s = t/t_{op}$, $d_{max}$ the target displacement and $t_{op}$ the operational time of the tracking strategy. The constant $d_{max}$ is the distance the fault slides dynamically in order to reach its sequent stable equilibrium point. Notice that the parameter $t_{op}$ is free to be decided depending on the earthquake control scenario that one wants to apply. Nevertheless, $t_{op}$ has to be sufficiently high with respect to the characteristic time of the earthquake phenomenon, but low enough to achieve aseismic slip with  higher velocity than the far-field velocity ($v_\infty$ in \eqref{eq:springmodel}), for the control scenario to make sense.

The choice of the reference output $y=x_1$ is motivated by the need to control the average slip over the fault. This average slip is directly connected with the magnitude of an earthquake through the seismic moment \cite{b:Kanamori-Brodsky-2004}. Therefore, by controlling the rate of the average slip, the system is forced to release its energy in a quasi-static way, \textit{i.e.}, aseismically. See \cite{b:Stefanou2019,b:https://doi.org/10.1029/2021JB023410} for more details.


\subsection{Sliding Mode based Control}
\label{ssec:sliding-modes}

To perform the tracking of the desired references $r(t),\dot{r}(t)$, a sliding-mode-based control is designed. Defining the tracking error variables,
\begin{equation}
  e_1=x_1-r, \quad e_2=x_2-\dot{r},
\end{equation}
the error dynamics reads as
\begin{equation}
\begin{split}
  \dot{e}_1&=e_2,\\
  \dot{e}_2&=-[\mu(e_1+r,e_2+\dot{r})-\mu^*]\hat{N}\sigma_n^\prime + \mu(e_1+r)\hat{N} p \\
   & \quad -\hat{k}(e_1+r) - \hat{\eta}(e_2+\dot{r}) + \hat{\varphi}_e(e_1+r,e_2+\dot{r},t) -\ddot{r} .
  \label{eq:esystem}
\end{split}  
\end{equation}

If the exact knowledge of the system parameters and the system dynamics would be available, all the known dynamics in $\dot{e}_2$ can be compensated, in order to get the nominal error system,
\begin{equation}
\begin{split}
  \dot{e}_1=e_2,\quad
  \dot{e}_2=\nu,
  \label{eq:esystem2}
\end{split} 
\end{equation}
by designing the control $p$ as
\begin{equation}
\begin{split}
  p &= \frac{1}{\mu(e_1+r,e_2+\dot{r})\hat{N}}\Big\{ \nu + [\mu(e_1+r)-\mu^*]\hat{N}\sigma_n^\prime \\ 
  & \quad + \hat{k}(e_1+r) + \hat{\eta}(e_2+\dot{r})  -\hat{\varphi}_e(e_1+r,e_2+\dot{r},t) + \ddot{r}\Big\},
  \label{eq:p}
\end{split}
\end{equation}
and with the new control input, $\nu$, designed to force $e_1,e_2$ towards zero. A solution could be the linear feedback control $\nu=-k_1 e_1 - k_2 e_2$ with any $k_1,k_2>0$. However, system \eqref{eq:esystem2} is valid only in the nominal case. If this is not the scenario, the application of the state-feedback \eqref{eq:p} with uncertain parameters will not lead to \eqref{eq:esystem2}. In this case, consider the feedback control
\begin{equation}
  p = \frac{1}{\mu_0 \hat{N}_0}\nu,
  \label{eq:ps}
\end{equation}
where the sub index `0' indicates the nominal value of the real considered parameter. Notice that some additional nominal parameters could be used in \eqref{eq:ps}, as $\hat{k}_0,\hat{\eta}_0$ or even the known term $\ddot{r}(t)$. However, the objective here is to design a controller requiring a limited amount of information.

Therefore, the closed-loop system obtained from \eqref{eq:esystem} and \eqref{eq:ps} reads as
\begin{equation}
\begin{split}
  \dot{e}_1=e_2,\quad
  \dot{e}_2=\beta(t,e)\left[\nu + h(t,e) \right],
  \label{eq:esystem3}
\end{split} 
\end{equation}
where $e=[e_1,e_2]^T$, $\beta(t,e)$, is the uncertain control coefficient, and $h(t,e)$ is a matched perturbation affecting the system. These terms read as
{\small
\begin{equation}
\begin{split}
  \beta(t,e) &= \frac{\mu(e_1+r,e_2+\dot{r})\hat{N}}{\mu_0 \hat{N}_0}, \\
  h(t,e) &= \frac{1}{\beta(t,e)}\Big\{-[\mu(e_1+r,e_2+\dot{r})-\mu^*]\hat{N}\sigma_n^\prime \\ & \quad -\hat{k}(e_1+r) - \hat{\eta}(e_2+\dot{r})+\hat{\varphi}_e(e_1+r,e_2+\dot{r},t)-\ddot{r} \Big\}.
\end{split}
\end{equation}}

Both these terms are assumed to fulfil in the operating domain
\begin{align}
0 < b_m & \leq \beta(t,e) \leq b_M\,, & \left| \frac{dh(t,e)}{dt}\right| & \leq \bar{L}\,,
\label{eq:Bounds}
\end{align} 
with known constants $b_m,b_M,\bar{L}$. 

\begin{remark}
The condition for $\beta(t,\,e)$ in \eqref{eq:Bounds} is satisfied because of the assumption of $\mu(e_1+r,e_2+\dot{r})=\mu(x_1,x_2)$ in \eqref{eq:mug}. The condition for $h(t,e)$ in \eqref{eq:Bounds} is satisfied (locally inside of a domain) because of the definition of $r(t)$ in \eqref{eq:ref} and if the external perturbation term $\hat{\varphi}_e(e_1+r,e_2+\dot{r},t)$ is Lipschitz w.r.t. time. As a result, the tracking result obtained in the sequel is valid locally.
\label{rem:domain}
\end{remark}

The design of the control input $\nu$ able to stabilize \eqref{eq:esystem3} at $e_1=e_2=0$, despite the presence of $\beta(t,\,e),h(t,e)$, results in an aseismic motion of system \eqref{eq:shift}. For this purpose, consider the Second-Order Continuous Twisting Algorithm (2-CTA) introduced in \cite{b:Torres-Sanchez-Fridman-Moreno}
\begin{equation}
	\begin{split}
	\nu = & - \lambda^{\frac{2}{3}} k_1\Sabs{e_1}^{\frac{1}{3}}- \lambda^{\frac{1}{2}} k_2\Sabs{e_2}^{\frac{1}{2}} + \xi_1,\\
	\dot{\xi}_1 = & - \lambda k_3\Sabs{e_1}^{0}- \lambda k_4\Sabs{e_2}^{0},
	\end{split}
	\label{eq:cta}
\end{equation}
and the Second-Order Discontinuous Integral Algorithm (2-DIA) introduced in \cite{b:Moreno_2016}
\begin{equation}
	\begin{split}
	\nu = & - \lambda^{\frac{1}{2}} k_{I2}\Sabs{\Sabs{e_2}^{\frac{3}{2}} + \lambda^{\frac{1}{2}} k_{I1}^{\frac{3}{2}} e_1}^{\frac{1}{3}} + \xi_1,\\
	\dot{\xi}_1 = & - \lambda k_{I3} \Sabs{e_1 + \lambda^{-\frac{1}{2}} k_{I4} \Sabs{e_2}^{\frac{3}{2}}}^{0}.
	\end{split}
	\label{eq:dia}
\end{equation} 

Both of these algorithms consist of a static homogeneous finite-time controller and a discontinuous integral action, aimed at estimating and compensating the uncertainties and perturbations. Notice that the presence of the discontinuous function, $\Sabs{\cdot}^{0}$, in the integral action finally results in a continuous control signal. 

\begin{theorem}
The origin of the system \eqref{eq:esystem3} is LFTS, despite the presence of the Lipschitz w.r.t. the time uncertainties/perturbations $h(t,e)$ and bounded uncertain coefficient $\beta(t,e)$ satisfying \eqref{eq:Bounds}, if the control $\nu$ takes the form of \eqref{eq:cta} or \eqref{eq:dia}, with gains appropriately chosen. 
\label{th:sliding}
\end{theorem}

As a consequence of Theorem \ref{th:sliding}, the state variables $\delta,\dot{\delta}$ of system \eqref{eq:springmodel}, are locally driven in finite-time to the desired references $r(t),\dot{r}(t)$ defined in \eqref{eq:ref}.

\begin{proof}
The closed-loop system \eqref{eq:esystem3}, with controller \eqref{eq:cta} reads as
\begin{equation}
\begin{split}
\dot{e}_1 &= e_2, \\
\dot{e}_2 &= \beta(t,e) \left(- \lambda^{\frac{2}{3}} k_1\Sabs{e_1}^{\frac{1}{3}}- \lambda^{\frac{1}{2}} k_2\Sabs{e_2}^{\frac{1}{2}} + e_{3}\right),\\
\dot{e}_{3} &= - \lambda k_3\Sabs{e_1}^{0}- \lambda k_4\Sabs{e_2}^{0} + \dot{h}(t,e) \,,
\label{eq:clsystem1}
\end{split}
\end{equation}
with $e_3 = \xi_1 + h(t,e)$. If the controller \eqref{eq:dia} is used, one gets
\begin{equation}
\begin{split}
\dot{e}_1 &= e_2, \\
\dot{e}_2 &= \beta(t,e) \left(- \lambda^{\frac{1}{2}} k_{I2}\Sabs{\Sabs{e_2}^{\frac{3}{2}} + \lambda^{\frac{1}{3}} k_{I1}^{\frac{3}{2}} e_1}^{\frac{1}{3}} + e_{3}\right),\\
\dot{e}_{3} &= - \lambda k_{I3} \Sabs{e_1 + \lambda^{-\frac{1}{2}} k_{I4} \Sabs{e_2}^{\frac{3}{2}}}^{0} + \dot{h}(t,e) \,.
\label{eq:clsystem2}
\end{split}
\end{equation}

The solutions of both systems are understood in the Filippov's sense (see \cite{b:filippov}). It will be shown that given bound parameters $b_{m}$, $b_{M}$ and a fixed $\bar{L}$, it is possible to find the values of both control gains such that $e=0$ is LFTS, but for a not assignable value $L^{*}$ of the Lipschitz bound $\bar{L}$. To meet this value, one just has to scale the previously obtained gains using $\lambda$ such that $\lambda L^{*}\geq \bar{L}$. 

The scaling with $\lambda>0$ does not alter the stability of the system. This can be shown by first performing the linear change of variables $z=\lambda e$ in systems \eqref{eq:clsystem1} and \eqref{eq:clsystem2}. The system in the variables $z$ has exactly the same form as \eqref{eq:clsystem1} and \eqref{eq:clsystem2}, respectively, but with the gains scaled as in \eqref{eq:cta}, \textit{i.e.}, $\left(k_{1},\,k_{2},\,k_{3},\,k_{4}\right)\rightarrow\left(k_{1}\lambda^{\frac{2}{3}},\,k_{2}\lambda^{\frac{1}{2}},\,k_{3}\lambda,\,k_{4}\lambda\right)$ and as in \eqref{eq:dia}, \textit{i.e.}, $\left(k_{I1},\,k_{I2},\,k_{I3},\,k_{I4}\right)\rightarrow\left(k_{I1}\lambda^{\frac{1}{3}},\,k_{I2}\lambda^{\frac{1}{2}},\,k_{I3}\lambda,\,k_{I4}\lambda^{-\frac{1}{2}}\right)$ (see \cite{b:Torres-Sanchez-Fridman-Moreno,b:Gutierrez-Mercado-Moreno-Fridman-IJRNC2020} for more details on each scaling). Therefore, both systems are equivalent in terms of stability.

The details of the proof for the closed-loop dynamics of system \eqref{eq:clsystem2} with the controller \eqref{eq:dia}, which follows closely the idea in \cite{b:MERCADOURIBE2020109018,b:Gutierrez-Mercado-Moreno-Fridman-IJRNC2020}, is presented. The proof for the controller \eqref{eq:cta} is similar, and a brief comment below is provided. First, the gain $\lambda=1$ is fixed. Then, a homogeneous and continuously differentiable Lyapunov function candidate is considered as
\begin{equation*}
V(\phi)=\frac{3}{5}\gamma_{1}\abs{\phi_{1}}^{\frac{5}{3}}+\phi_{1}\phi_{2}+\frac{2}{5}k_{I1}^{-\frac{3}{2}}\abs{\phi_{2}}^{\frac{5}{2}}+\frac{1}{5}\abs{\phi_{3}}^{5}\,, 
\end{equation*}
where 
\begin{equation*}
\phi_{1}= e_{1}-\Sabs{\phi_{3}}^{3}\,,\quad 
\phi_{2}=  e_{2}\,,\quad
\phi_{3}=  k_{I1}^{-\frac{1}{2}}k_{I2}^{-1}e_{3}\,.
\end{equation*}
$V\left(\phi\right)$ is positive definite if $k_{I1}>0$ and $\gamma_{1}>0$ is selected sufficiently large. Its derivative along the trajectories of the system \eqref{eq:clsystem2} is
{\small
\begin{align*}
\dot{V} & =-k_{I1}^{\frac{1}{2}}k_{I2}\beta(t,e)\left(\phi_{1}+k_{I1}^{-\frac{3}{2}}\Sabs{\phi_{2}}^{\frac{3}{2}}\right)\\ 
&\quad \times \left[\Sabs{\left(\phi_{1}+k_{I1}^{-\frac{3}{2}}\Sabs{\phi_{2}}^{\frac{3}{2}}\right)+\Sabs{\phi_{3}}^{3}}^{\frac{1}{3}}-\phi_{3}\right]\\
 &\quad + \left(\gamma_{1}\Sabs{\phi_{1}}^{\frac{2}{3}}+\phi_{2}\right)\phi_{2}-k_{I3}k_{I1}^{-\frac{1}{2}}k_{I2}^{-1}\Big[\Sabs{\phi_{3}}^{4} \\
 &\quad -3\abs{\phi_{3}}^{2}\left(\gamma_{1}\Sabs{\phi_{1}}^{\frac{2}{3}}+\phi_{2}\right)\Big]\Bigg\{ \Sabs{\phi_{1}+\Sabs{\phi_{3}}^{3}+k_{I4}\Sabs{\phi_{2}}^{\frac{3}{2}}}^{0}\\
 &\quad -\frac{1}{k_{I3}}\dot{h}(t,e)\Bigg\} \,.
\end{align*}}
Using the bounds for $\beta(t,e)$ and $h(t,e)$ in \eqref{eq:Bounds}, the derivative reads as
{\small
\begin{align*}
\dot{V} &\leq-\tilde{k}_{I2}b_{m}\left(\phi_{1}+k_{I1}^{-\frac{3}{2}}\Sabs{\phi_{2}}^{\frac{3}{2}}\right)\\
&\quad \times \left[\Sabs{\left(\phi_{1}+k_{I1}^{-\frac{3}{2}}\Sabs{\phi_{2}}^{\frac{3}{2}}\right)+\Sabs{\phi_{3}}^{3}}^{\frac{1}{3}}-\phi_{3}\right]\\
 &\quad + \left(\gamma_{1}\Sabs{\phi_{1}}^{\frac{2}{3}}+\phi_{2}\right)\phi_{2}-\tilde{k}_{I3}\Big[\Sabs{\phi_{3}}^{4} \\
 &\quad -3\abs{\phi_{3}}^{2}\left(\gamma_{1}\Sabs{\phi_{1}}^{\frac{2}{3}}+\phi_{2}\right)\Big]\Bigg\{ \Sabs{\phi_{1}+\Sabs{\phi_{3}}^{3}+k_{I4}\Sabs{\phi_{2}}^{\frac{3}{2}}}^{0}\\
 &\quad -\left[-\tilde{L}^{*},\,\tilde{L}^{*}\right]\Bigg\} \,,
\end{align*}}
where $\tilde{L}^{*}=\frac{L^{*}}{k_{I3}}$, $\tilde{k}_{I2}=k_{I1}^{\frac{1}{2}}k_{I2}$
and $\tilde{k}_{I3}=k_{I3}k_{I1}^{-\frac{1}{2}}k_{I2}^{-1}$. Note that the value of $L^{*}$ is not given, but it has to be found. Using Lemma \ref{lem:Holder S}, the first term can be proven to be negative semi-definite and it vanishes only on the set $S_{1}=\conjunto{k_{I1}^{\frac{3}{2}}\phi_{1}+\Sabs{\phi_{2}}^{\frac{3}{2}}=0}$.

Evaluating $\dot{V}$ on this set results in
{\small
\begin{align*}
\Subsin{\dot{V}}_{S_{1}}  &\leq-\left(\gamma_{1}k_{I1}^{-1}-1\right)\left|\phi_{2}\right|^{2}\\
&\quad -\tilde{k}_{I3}\abs{\phi_{3}}^{2}\left[3\left(\gamma_{1}k_{I1}^{-1}-1\right)\phi_{2}+\Sabs{\phi_{3}}^{2}\right] \\
&\quad \times \left\{ \Sabs{\left(-k_{I1}^{-\frac{3}{2}}+k_{I4}\right)\Sabs{\phi_{2}}^{\frac{3}{2}}+\Sabs{\phi_{3}}^{3}}^{0}-\left[-\tilde{L}^{*},\,\tilde{L}^{*}\right]\right\} \,.
\end{align*}}

If $\gamma_{1}>k_{I1}$, the first term is negative semidefinite, and it is zero only on the set $S_{2}=\conjunto{\phi_{2}=0}$. Evaluating $\Subsin{\dot{V}}_{S_{1}}$ on $S_{2}$ reads as 
\begin{align*}
\Subsin{\dot{V}}_{S_{1}\cap S_{2}} & \leq-\tilde{k}_{I3}\Sabs{\phi_{3}}^{4}\left\{ \Sabs{\phi_{3}}^{0}-\left[-\tilde{L}^{*},\,\tilde{L}^{*}\right]\right\} \,.
\end{align*}
This is negative if $\tilde{L}^{*}=\frac{L^{*}}{k_{I3}}<1$. By Lemma \ref{lem:Negative I Semi}, $\Subsin{\dot{V}}_{S_{1}}<0$ is proven by selecting $\tilde{k}_{I3}>0$ small. Using Lemma \ref{lem:Negative I Semi}, again, it is possible to make $\dot{V}<0$ selecting $\tilde{k}_{I2}>0$ sufficiently large.

For the closed-loop system \eqref{eq:clsystem1}, the smooth and homogeneous Lyapunov function
\begin{align*}
V(e) &= \alpha_{1}\abs{e_{1}}^{\frac{5}{3}}+\alpha_{2}e_{1}e_{2}+\alpha_{3}\abs{e_{2}}^{\frac{5}{2}}+\alpha_{4}e_{1}\Sabs{e_{3}}^{2}\\
&\quad -\alpha_{5}e_{2}e_{3}^{3}+\alpha_{6}\abs{e_{3}}^{5},
\end{align*}
can be selected. It has been shown in \cite{b:Torres-Sanchez-Fridman-Moreno} that the origin of system \eqref{eq:clsystem1} is GFTS for appropriate selected gains $k_{1}>0,\,k_{2}>0,\,k_{3}>0,\,k_{4}$. Although this result has been obtained in \cite{b:Torres-Sanchez-Fridman-Moreno} using an SOS algorithm and assuming that the control coefficient is known and constant, the proof can be also extended to the actual case using similar arguments as those used above for the controller \eqref{eq:dia}. 

Therefore, the origin of systems \eqref{eq:clsystem1} and \eqref{eq:clsystem2} is asymptotically stable, but locally, as a consequence of Remark \ref{rem:domain}. Moreover, since the systems are homogeneous of negative degree, such origins are LFTS (see \cite{b:levant_automatica2005}).

\end{proof}

\subsection{LQR-based Control}

Following the sliding-mode based control design, this subsection presents an extended Linear Quadratic Regulator (e-LQR). The term \textit{extended} is due to the integral action added to a standard LQR algorithm, in such a way that Lipschitz w.r.t. the states uncertainties are compensated with the resultant control. 

Starting from \eqref{eq:shift}, the plant is extended with a double integrator to improve the tracking of the reference trajectory. The stability of this augmented plant controlled by a LQR controller is being proved by a Lyapunov's approach. For this purpose, define
\begin{equation}
\begin{split}
  \dot{\xi}_1 = x_1 - r, \quad
  \dot{\xi}_2 = \xi_1 .
\end{split}
\label{eq:double_integrator}
\end{equation}

The matrix form of system \eqref{eq:shift} is
\begin{equation}
\begin{split}
\displaystyle
\underbrace{\begin{bmatrix}
\dot{x}_1 \\
\dot{x}_2
\end{bmatrix}}_{\dot{x}} &= 
\underbrace{\begin{bmatrix}
0 & 1 \\
-\hat{k} & -\hat{\eta}
\end{bmatrix}}_{A(t)}
\underbrace{\begin{bmatrix}
x_1 \\
x_2
\end{bmatrix}}_{x} + 
\underbrace{\begin{bmatrix}
0 \\
\mu(x_1,x_2) \hat{N}
\end{bmatrix}}_{B(t,x)} 
p \\ &\quad +
\underbrace{\begin{bmatrix}
0 \\
- [\mu(x_1,x_2) - \mu^* ] \hat{N} \sigma^{\, \prime}_{\mathrm{n}}+\hat{\varphi}_e(x_1,x_2,t)
\end{bmatrix},}_{g(t,x)}
\end{split}
\label{eq:system_mf_lqr}
\end{equation}
whereas the matrix form of $(\xi_1,\xi_2)$-system is
\begin{equation}
\displaystyle
\underbrace{\begin{bmatrix}
\dot{\xi}_1 \\
\dot{\xi}_2
\end{bmatrix}}_{\dot{\xi}} = 
\underbrace{\begin{bmatrix}
0 & 0 \\
1 & 0
\end{bmatrix}}_{C_{\xi}}
\underbrace{\begin{bmatrix}
\xi_1 \\
\xi_2
\end{bmatrix}}_{\xi} +
\underbrace{\begin{bmatrix}
1 & 0 \\
0 & 0
\end{bmatrix}}_{C_{x}}
\underbrace{\begin{bmatrix}
x_1 \\
x_2
\end{bmatrix}}_{x} 
+
\underbrace{\begin{bmatrix}
-r(t) \\
0
\end{bmatrix}}_{r_\xi(t)} .
\label{eq:double_integrator_mf}
\end{equation}

Following a conventional integral control design (see \cite[Chapter 12]{b:Khalil2002} for example) and supposing a constant reference $r(t)=r_0$, the augmented system composed by \eqref{eq:system_mf_lqr}-\eqref{eq:double_integrator_mf} reads as
\begin{equation}
\displaystyle
\underbrace{\begin{bmatrix}
\dot{x} \\
\dot{\xi}
\end{bmatrix}}_{\dot{x}_a} = 
\underbrace{\begin{bmatrix}
A(t) &  0_{2\times 2} \\
C_x & C_\xi
\end{bmatrix}}_{A_a(t)}
\underbrace{\begin{bmatrix}
x \\
\xi
\end{bmatrix}}_{x_a} + 
\underbrace{\begin{bmatrix}
B(t,x) \\
0_{2\times 1}
\end{bmatrix}}_{B_a(t,x_a)} 
p +
\underbrace{\begin{bmatrix}
g(t,x) \\
0_{2\times 1} 
\end{bmatrix}}_{g_a(t,x_a)}.
\label{eq:system_augm_lqr}
\end{equation}

Consider additive (matched) uncertainties to $A_a(t)$ and $B_a(t,x_a)$ such that,
\begin{align}\label{eq:DA}
\displaystyle
A_a(t) &= A_0 + \Delta A(t), \\
B_a(t,x_a) &= B_0 + \Delta B(t,x_a),
\label{eq:DB}
\end{align} 
where
{\small
\begin{align*}
\displaystyle
&A_0 =
\begin{bmatrix}
0 & 1 & 0 & 0 \\
-\hat{k}_0 & -\hat{\eta}_0 & 0 & 0 \\
1 & 0 & 0 & 0 \\
0 & 0 & 1 & 0
\end{bmatrix}, 
\Delta B (t,x_a) = \begin{bmatrix}
0 \\
\Delta [\mu(x_1,x_2) \hat{N}] \\
0 \\
0
\end{bmatrix} \\
&\Delta A(t) =
\begin{bmatrix}
0 & 0 & 0 & 0 \\
-\Delta \hat{k} & -\Delta \hat{\eta} & 0 & 0 \\
0 & 0 & 0 & 0 \\
0 & 0 & 0 & 0
\end{bmatrix}, 
B_0 = \begin{bmatrix}
0 \\
[\mu_{res} \hat{N} \\
0 \\
0
\end{bmatrix}.
\end{align*}}

The sub-index `0' represents the nominal value, whereas the quantity with the prefix `$\Delta$' corresponds to the uncertainties of the respective variable. In addition, the assumption \eqref{eq:mug} was used for the friction coefficient $\mu(x_1,x_2)$, so the variation $\Delta B (t,x_a)$ is always positive semi-definite.

From \eqref{eq:system_augm_lqr}, \eqref{eq:DA}, and \eqref{eq:DB}, one gets
\begin{equation}
\displaystyle
\dot{x}_a = A_0x_a + \Delta B(t,x_a)p + B_0 \left[p + h(t,x_a)\right],
\label{eq:system_mu_lqr}
\end{equation}
where $\displaystyle h(t,x_a) = B_0^+ \Delta A(t)x_a + B_0^+ g(t,x_a)$, with $\displaystyle B_0^+$ is the \textit{Moore-Penrose inverse matrix} of $\displaystyle B_0$, \textit{i.e.},
\begin{equation}
\displaystyle
B\strut^+_0=\begin{bmatrix}
0 & \frac{1}{\mu_{res} \hat{N}} & 0 & 0
\end{bmatrix}.
\label{eq:B0+}
\end{equation}

Then, the nonlinear vector $h(t,x_a)$ can be written as
\begin{equation}
\begin{split}
h(t,x_a) &= B\strut^+_0 \left[\Delta A(t) x_a + g(t,x_a) \right] \\
&= -\frac{1}{\mu_{res} \hat{N}} \Bigg\{\Delta \hat{\eta} x_2 + \Delta \hat{k} x_1 \\ &\quad + \left[ \mu(x_1,x_2) - \mu^* \right]  \sigma'_{\mathrm{n}} \hat{N} -\hat{\varphi}_e(x_1,x_2,t)\Bigg\}.
\end{split}
\end{equation}

Assuming the external perturbation $\hat{\varphi}_e(x_1,x_2,t)$ to be Lipschitz w.r.t. the states, \textit{i.e.}, $\abs{\hat{\varphi}_e(x_1,x_2,t)}\leq \abs{\hat{\varphi}_{1e}x_1 + \hat{\varphi}_{2e}x_2}$ for some known positive constants $\hat{\varphi}_{1e},\hat{\varphi}_{2e}$, the norm of $h(t,x_a)$ reads as
{\small
\begin{align*}
\displaystyle
&\left\Vert h(t,x_a) \right\Vert  = \\
&= \Bigg\| \frac{1}{\mu_{res} \hat{N}} \\ & \quad \times
\begin{bmatrix}
\Delta \hat{k}  + \frac{\mu(x_1,x_2)  - \mu^*}{x_1}  \sigma'_{\mathrm{n}} \hat{N} + \hat{\varphi}_{1e} & \Delta \hat{\eta} + \hat{\varphi}_{2e} & 0 & 0
\end{bmatrix}
x_a
\Bigg\| 
\end{align*}}
\vspace{-5pt}
{\footnotesize
\begin{align*}
\displaystyle
&\leq \Bigg\| \frac{1}{\mu_{res} \hat{N}} \\ & \quad \times
\begin{bmatrix}\Delta \hat{k}_{\mathrm{max}} + \mu^{\mathrm{max}}_{x_1} \sigma'_{\mathrm{n}} \Delta \hat{N}_{\mathrm{max}} + \hat{\varphi}_{1e} & \Delta \hat{\eta}_{\mathrm{max}} + \hat{\varphi}_{2e} & 0 & 0 \end{bmatrix} x_a \Bigg\|,
\end{align*}}
\hspace{-5pt}where the bound $\mu^{\mathrm{max}}_{x_1}=\left| \frac{\partial \mu(x_1,x_2) }{\partial x_1} \right|_{\mathrm{max}}$ corresponds to the maximum absolute softening slope of the friction. The subscript `max' denotes the maximum variation in absolute term from the respective nominal values.

Finally, one gets the bounds of the variation coefficient $\Delta B(t,x_a)$ and perturbation $h(t,x_a)$ in system \eqref{eq:system_mu_lqr} as,
\begin{align}
0 &\leq \Delta B(t,x_a) \,, & \left\Vert h(t,x_a) \right \Vert &\leq \left\Vert G x_a \right \Vert \,,
\label{eq:lipschitz1}
\end{align}
with $G$ defined as,
{\small
\begin{equation}
\begin{split}
G &= \frac{1}{[\mu(x_1) \hat{N}]_{min}} \\ 
& \quad \times \begin{bmatrix}\Delta \hat{k}_{\mathrm{max}} + \mu^{\mathrm{max}}_{\tilde{x}_1} \sigma'_{\mathrm{n}} \Delta \hat{N}_{\mathrm{max}} + \hat{\varphi}_{1e} & \Delta \hat{\eta}_{\mathrm{max}} + \hat{\varphi}_{2e} & 0 & 0 \end{bmatrix}.
\end{split}
\label{eq:lipschitz2}
\end{equation}}   

Inspired from the original LQR control proposed in \cite{b:Kalman-1960}, define the e-LQR control input $p$ designed for the augmented system \eqref{eq:system_mu_lqr} as
\begin{equation}
\displaystyle
p = -R^{-1}B_0^{\text{T}} \Theta x_a = -[k_1 \quad k_2 \quad k_3 \quad k_4][x_1 \quad x_2 \quad \xi_1 \quad \xi_2]^T,
\label{eq:u}
\end{equation}
where $R$ is a \textit{positive definite} matrix to be chosen and $\Theta$ the \textit{positive-definite} solution of the following Continuous Algebraic Riccati Equation (CARE),
\begin{equation}
\displaystyle
A_0^{\text{T}} \Theta + \Theta A_0 - \Theta B_0 R^{-1} B_0^{\text{T}} \Theta = -Q,
\label{eq:CARE}
\end{equation}
with
\begin{equation}
\displaystyle
Q = Q_0 + ||R^{1/2}||\strut^2 G^{\text{T}} G ,
\label{eq:Q}
\end{equation}
and $Q_0$ a \textit{positive definite} matrix to be chosen. 

The next Theorem states the main result of the e-LQR control.

\begin{theorem}
The origin of the augmented closed-loop system \eqref{eq:system_mu_lqr} is GES, with the presence of the Lipschitz w.r.t. the states perturbation $h(t,x_a)$ and positive semi-definite term $\Delta B(t,x_a)$ satisfying \eqref{eq:lipschitz1} and \eqref{eq:lipschitz2}, when the control input $p$ takes the form of \eqref{eq:u}-\eqref{eq:Q}.
\label{th:eLQR}
\end{theorem}

As a consequence of Theorem \ref{th:eLQR}, the slip $\delta$ and slip-rate fault $\dot{\delta}$ of the original system \eqref{eq:springmodel} are driven globally and exponentially to a desired constant reference $r(t)=r_0$.

\begin{proof}
Select as Lyapunov candidate the positive definite and radially unbounded function $\displaystyle V(x_a)=x_a^{\mathrm{T}} \Theta x_a$, where $\displaystyle \Theta=\Theta^{\mathrm{T}}>0_{4\times 4}$ as a positive definite matrix. The following inequalities are true for the Lyapunov candidate
\begin{align}
\displaystyle
\lambda_{\text{min}}(\Theta)||x_a||^2 \leq V(x_a) \leq \lambda_{\text{max}}(\Theta)||x_a||^2, \\
\left|\left| \frac{\partial V(x_a)}{\partial x_a} \right| \right| \leq 2 \lambda_{\text{max}}(\Theta) ||x_a||,
\label{eq:lyap_ineq}
\end{align}
where $\lambda_{\text{min}}(\Theta)$ and $\lambda_{\text{max}}(\Theta)$ are the minimum and maximum eigenvalues of the matrix $\Theta$, respectively. 

The time derivative of the Lyapunov candidate $\displaystyle V(x_a)$ along the trajectories of closed loop system \eqref{eq:system_mu_lqr} and \eqref{eq:u} reads as,
\begin{subequations}
\begin{align}
\displaystyle
\dot{V} &=\dot{x_a}^{\mathrm{T}} \Theta x_a + x_a^{\mathrm{T}} \Theta \dot{x_a} \nonumber \\
&= x_a^{\mathrm{T}} \left(A^{\mathrm{T}}_0 \Theta + \Theta A_0 -\Theta B_0 R^{-1} B^{\mathrm{T}}_0 \Theta \right) x_a  \label{eq:lyap_deriva} \\
& \quad - 2 x_a^{\mathrm{T}} \Theta \frac{B_0 R^{-1} \Delta B^{\mathrm{T}} + \left(B_0 R^{-1} \Delta B^{\mathrm{T}}\right)^{\mathrm{T}}}{2} \Theta x_a  \label{eq:lyap_derivb} \\
& \quad -x_a^{\mathrm{T}}\Theta B_0 R^{-1} B^{\mathrm{T}}_0 \Theta x_a + 2x_a^{\mathrm{T}} \Theta B_0 h,
\end{align}
\label{eq:lyap_deriv}\end{subequations}
where the dependencies of $h(t,x_a),\Delta B(t,x_a)$ have been omitted for simpler notation.

The term (\ref{eq:lyap_deriva}) is the CARE as defined in \eqref{eq:CARE}. Moreover, due to the fact that the nominal matrix $\displaystyle B_0$ has been selected to have the term $\Delta B(t,x_a)$ non-negative, the non-Hermitian matrix $\displaystyle B_0 R^{-1} \Delta B^{\mathrm{T}}$ is always \textit{positive semi-definite} and the term on (\ref{eq:lyap_derivb}) is a non-positive scalar. 

Therefore, letting $\displaystyle z=R^{-1/2} B^{\mathrm{T}}_0 \Theta x_a$, one gets
\begin{equation}
\displaystyle
\dot{V} \leq - x_a^{\mathrm{T}} Q x_a - z^{\mathrm{T}}  z + 2 z^{\mathrm{T}} R^{1/2} h,
\label{eq:lyap_deriv2}
\end{equation}
with $\displaystyle R^{-1}=R^{-1/2}R^{-1/2}$. 

Using the definition of $Q$ in (\ref{eq:Q}) and the Lipschitz condition for $h$ in \eqref{eq:lipschitz1}, it comes
\begin{equation}
\begin{split}
\dot{V} \leq &- x_a^{\mathrm{T}} Q_0 x_a - \left\Vert R^{1/2} \right\Vert^2 x_a^{\mathrm{T}} G^{\text{T}} Gx_a - z^{\mathrm{T}}  z + 2 z^{\mathrm{T}} R^{1/2} h, \\
\leq &- x_a^{\mathrm{T}} Q_0 x_a - \left\Vert R^{1/2} \right\Vert^2 \left\Vert Gx_a \right\Vert^2 - \left\Vert z \right\Vert^2 \\ & + 2 \left\Vert R^{1/2} \right\Vert \left\Vert z \right\Vert \left\Vert Gx_a \right\Vert, \\
\leq &- x_a^{\mathrm{T}} Q_0 x_a - \left( \left\Vert R^{1/2} \right\Vert \left\Vert Gx_a \right\Vert - \left\Vert z \right\Vert \right)^2, \\
\leq &- x_a^{\mathrm{T}} Q_0 x_a \leq -\lambda_{\text{min}}(Q_0) ||x_a||^2 < 0.
\end{split}
\end{equation}

As a conclusion, $V(x_a)$ is a Lyapunov function for system (\ref{eq:system_mu_lqr}). Therefore, its origin is GES. 
\label{pr:eLQR}
\end{proof}

\subsection{Control Strategies Comparison}

The designed control strategies, $p$, shown in \eqref{eq:ps}, \eqref{eq:cta}, \eqref{eq:dia} and \eqref{eq:u}-\eqref{eq:Q} represent the fluid injected to the fault able to achieve aseismic response in the model \eqref{eq:springmodel} by tracking a slow reference, robustly and by using a continuous control signal. This is performed despite the presence of uncertainties and/or disturbances assumed to be as \eqref{eq:esystem3},\eqref{eq:Bounds} for the sliding mode algorithms, and to be as \eqref{eq:system_mu_lqr},\eqref{eq:lipschitz1},\eqref{eq:lipschitz2} for the e-LQR. For the earthquake application, this basically means that the friction coefficient must be Lipschitz w.r.t. the time and the states of the dynamical system. Such assumption is fulfilled by the most common friction laws used in fault mechanics (see \cite{b:https://doi.org/10.1029/2021JB023410} for a mathematical proof).

Now, some properties of each control strategy are discussed in the following, highlighting the differences between them.

Sliding mode-based control:
\begin{itemize}
  \item The origin of system \eqref{eq:esystem3}, with bounded control coefficient $\beta(t,e)$ and Lipschitz w.r.t. the time perturbation $h(t,e)$ assumed as \eqref{eq:Bounds}, is LFTS.
  \item Calculation of the gains for the 2-CTA \eqref{eq:cta} and 2-DIA \eqref{eq:dia} controllers, can be obtained using a Sum Of Square algorithm for the 2-CTA (see \cite{b:Torres-Sanchez-Fridman-Moreno}) and by performing a maximization of homogeneous functions for the 2-DIA (see \cite{b:MERCADOURIBE2020109018,b:Gutierrez-Mercado-Moreno-Fridman-IJRNC2020}).
  \item The reference signal $r(t)$ that can be tracked, has to be chosen as $\abs{{r}^{(3)}(t)}\leq \gamma$ with a positive constant $\gamma$, in order to fulfil the assumptions in \eqref{eq:Bounds}. Note that the selected reference \eqref{eq:ref} fulfils this condition for all $t\in [0,t_{op}]$.
  \item Systems \eqref{eq:clsystem1} and \eqref{eq:clsystem2} are homogeneous vector-set of degree $d=-1$ and weights $({r}_{1},{r}_{2},{r}_{3})=(3,2,1)$. Due to homogeneity properties \cite{b:levant_automatica2005}, the theoretical precision of the states after the transient are $\left| { e }_{ 1 } \right| < \Delta_1 T_s^3$, $\left| { e }_{ 2 } \right| < \Delta_2 T_s^2$ and $\left| { e }_{ 3 } \right| < \Delta_3 T_s$, where $\Delta_i>0$ with $i=\left\{1,...,3\right\}$ and $T_s$ the sampling time.
\end{itemize}

LQR-based control:
\begin{itemize}
  \item The origin of system \eqref{eq:system_mu_lqr}, with positive semi-definite variation coefficient $\Delta B(t,x_a)$ and Lipschitz w.r.t. the states $h(t,x_a)$ assumed as \eqref{eq:lipschitz1}, is GES.
  \item Calculation of the gains for the e-LQR control \eqref{eq:u} are obtained by solving the CARE \eqref{eq:CARE}-\eqref{eq:Q}.
  \item  The classical version of an integral control will track a constant reference, \textit{i.e.}, $r(t)=r_0$ (see \cite[Chapter 12]{b:Khalil2002}, for example). According to the internal model principle (see \cite{b:Francis-Wonhamm-1976}), the use of a double integrator \eqref{eq:double_integrator} will be able to follow linear time references, \textit{i.e.}, $r(t)=\alpha_1 t+\alpha_2$, with some $\alpha_1,\alpha_2 \in \mathbb{R}$. The steady-state error tracking of the target reference \eqref{eq:ref} using the presented control will not be zero, but it will be improved by using this double integrator scheme. In addition, this error can become smaller increasing the e-LQR integral gains.
\end{itemize}

\section{Numerical Simulations and Experimental Confirmation}
\label{sec:SimExp}

\begin{figure}[ht!]
  \centering 
  \includegraphics[width=8.5cm,height=4cm]{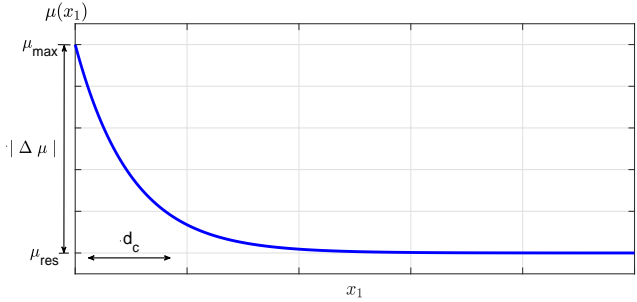}
  \caption{Slip-weakening friction law used in the simulations.}
  \label{fig:mu}
\end{figure}

In order to illustrate the performance of both previous control algorithms, simulations and laboratory experiments have been performed based on the shifted system described by \eqref{eq:shift}. The definition of the next parameters to be used are recalled
\begin{equation*}
\begin{split}
  \hat{N}&=\frac{A}{m}, \quad \hat{k}=\frac{k}{m}, \quad \hat{\eta}=\frac{\eta}{m}, \quad A=L_{ac}^2, \\
  k &= \frac{G}{L_{ac}}, \quad m = \rho L_{ac}^3, \quad \hat{\varphi}_e(x_1,x_2,t)=\frac{\varphi_e(x_1,x_2,t)}{m},
\end{split}
\end{equation*}
where $\varphi_e(x_1,x_2,t)$ is an external perturbation and $\mu(x_1,x_2)$ the friction coefficient. In the following numerical examples, a friction coefficient of the form 
\begin{equation}
  \mu(x_1)=\mu_{res}-\Delta \mu \cdot e^{-\nicefrac{x_1}{d_c}},
  \label{eq:mu}
\end{equation}
is considered with $\Delta \mu<0$. Such function is defined as a slip-weakening friction law \cite{b:Kanamori-Brodsky-2004} and it evolves from an initial value $\mu_{max}=\mu_{res}-\Delta \mu$ (static friction coefficient), to a residual one $\mu_{res}$ (kinetic friction coefficient) in a characteristic slip $d_c$, as shown by Fig. \ref{fig:mu}.

The choice of a friction law that depends only on the slip is argued by the difficulty of modelling the complex structure region of frictional interfaces. Nevertheless, the friction coefficient \eqref{eq:mu} fulfils the assumption imposed in \eqref{eq:mug} and it was used for calculating the controller gains.

Two cases have been studied in this paper. For the numerical simulations, the parameters of the reduced model for earthquakes (see eq. \eqref{eq:shift}) are chosen to correspond to an earthquake of magnitude $M_w=5.8$, where simulations are being made using MATLAB Simulink. For the experiments, a different set of parameters was chosen corresponding to the laboratory-fault setup (Fig. \ref{fig:spring-slider-lab}), corresponding to an earthquake of magnitude $M_w=-5.2$. The mechanical and frictional parameters for these two cases are given in Table \ref{tab:para}.
 
\begin{table*}[ht!]
\caption{Mechanical and frictional properties adopted for the real-scale simulations (R) and the laboratory experiments (L)}
\begin{center}
\begin{threeparttable}  
\begin{tabular}{|c|c|c|c|c|}
\hline 
\textbf{Parameter} & \textbf{Description} & \textbf{Simulations} & \textbf{Experiments}\tnote{1} & \textbf{Scaling factor}\tnote{2} \\ 
 &  & \textbf{Real-fault (R)} & \textbf{Lab-fault (L)} & \textbf{(R/L)} \\ 
\hline 
$\rho$ & density & 2500 [kg/m$^3$] & 1385 [kg/m$^3$] & $\lambda_\rho=1.81$ \\ 

$G$ & shear modulus & 30 [GPa] & 225.5 [kPa] & $\lambda_G=1.33\times 10^5$  \\ 
 
$\eta$ & damping coefficient & $5\times 10^{14}$ [kg/s] & 408 [kg/s] & $\lambda_\eta=1.23\times 10^{12}$ \\ 

$L_{ac}$ & activated fault length & 5 [km] & 0.1 [m] & $\lambda_{L_{ac}}=5\times 10^4$ \\ 

$\sigma_n^\prime$ & effective normal stress & 50 [MPa] & 0.1 [MPa] & $\lambda_p=500$ \\ 
 
$\mu_{res}$ & residual friction & 0.2353 & 0.4 & $\lambda_\mu=0.5882$ \\ 

$\Delta \mu$ & friction drop & -0.1 & -0.17 & $\lambda_\mu$ \\ 
 
$d_c$ & characteristic slip distance & 276.35 [mm] & 2.5 [mm] & $\lambda_\delta=110.54$ \\ 

$d_{max}$ & maximum displacement & 785 [mm] & 7.1 [mm] & $\lambda_\delta$ \\ 
 
$t_{op}$ & operation time & 184.17 [h] & 1 [h] & $\lambda_t=184.17$ \\ 

$T_s$ & sampling time & 184 [ms] & 1 [ms] & $\lambda_t$ \\ 
  
$M_0$ & seismic moment & $6.25\times 10^{17}$ [Nm] & 17 [Nm] & $\lambda_{M_0}=3.68\times 10^{16}$ \\ 
 
$M_w$ & seismic magnitude & 5.8 & -5.2 & $\lambda_{M_w}=4.97$ \\
\hline   
\end{tabular}
\begin{tablenotes}
\item[1] These values correspond to a double interface of an area of $A =100$ [cm$^2$] and a spring coefficient of $k=45.1$ [N/mm] used for the experimental apparatus depicted in Fig. \ref{fig:spring-slider-lab}.\\
\item[2] Due to the fact that seismic magnitude is based on a logarithm scale (see eq. \eqref{eq:magnitude}), the scaling factor between both faults is equal to $\lambda_{M_w}=\frac{2}{3} \log_{10}\lambda_{M_0}-6.07$. 
\end{tablenotes}
\end{threeparttable}
\end{center}
\label{tab:para}
\end{table*}

\begin{remark}
Based on the Buckingham $\pi$ theorem \cite{b:Logan-2013}, the lab-fault dynamics can be upscaled to obtain the real-fault dynamics under appropriate scaling laws (see \cite[Appendix E]{b:Tzortzopoulos-2021} for further details).
\end{remark}

Note how the seismic magnitude of the lab-fault allows to reproduce safely in the laboratory earthquake-like instabilities, but such experimental results can be upscaled to real earthquake events.
 
\begin{remark}
The controller gains used for the numerical simulations were obtained by upscaling the gains chosen for the experimental tests. This was made in order to be consistent with the upscale process commented above. See Appendix \ref{sec:app1} for more details.
\end{remark}

\subsection{Numerical Simulations}

The numerical simulations on the shifted system \eqref{eq:shift} using the real-fault parameters from Table \ref{tab:para}, have been performed in MATLAB Simulink with Dormand-Prince's integration method of a fixed time step equal to $T_s^R=184$ [ms].  

The controller \eqref{eq:ps} with $\mu_0=\mu_{res}^L, \hat{N}_0=A^L/m^L$, $\nu$ chosen as \eqref{eq:cta} for the 2-CTA and as \eqref{eq:dia} for the 2-DIA, and the e-LQR control \eqref{eq:u} have been implemented.\footnote{Superscript $L$ means lab-fault parameters from Table \ref{tab:para}.} The controllers gains have been selected as 
\begin{itemize}
  \item 2-CTA: $k_1=781.37$, $k_2=3.22\times 10^3$, $k_3=4.51 \times 10^{-4}$, $k_4=2.15 \times 10^{-4}$ and $\lambda=500$.
  \item 2-DIA: $k_{I1}=5.21\times 10^{-2}$, $k_{I2}=3.23\times 10^{3}$, $k_{I3}=3.91 \times 10^{-4}$, $k_{I4}=0$ and $\lambda=500$.
  \item e-LQR: $k_1=1.88\times 10^{9}$, $k_2=5.79\times 10^8$, $k_3=1.02\times 10^{6}$ and $k_4=18.52$.
  \label{simgains}
\end{itemize}

The reference signal \eqref{eq:ref} is used with $d_{max}=785$ [mm] and $t_{op}=184.17$ [h]. Furthermore, an external perturbation $\hat{\varphi}_e(x_1,x_2,t)= 3.2\times 10^{-4}\sin(0.69t) + 3.2\times 10^{-6}x_1$ [m/s$^2$] has been added to the system, and the initial condition $x_1(0)=x_2(0)=0$ was chosen for the three algorithms. It is important to notice that the external perturbation has a Lipschitz w.r.t. the time term and a Lipschitz w.r.t. the states term. Furthermore, the chosen frequency in the sinusoidal perturbation is equal to the natural frequency of the system, \textit{i.e.}, $\omega_n=\sqrt{\nicefrac{k}{m}}\simeq 0.69$ [rad/s]. The magnitude of both terms in such perturbation has been selected to appreciate a change in the system response and to check the robustness of the controllers versus different kinds of functions and resonance behaviour.

\begin{figure}[ht!]
  \centering 
  \includegraphics[width=8.3cm,height=4cm]{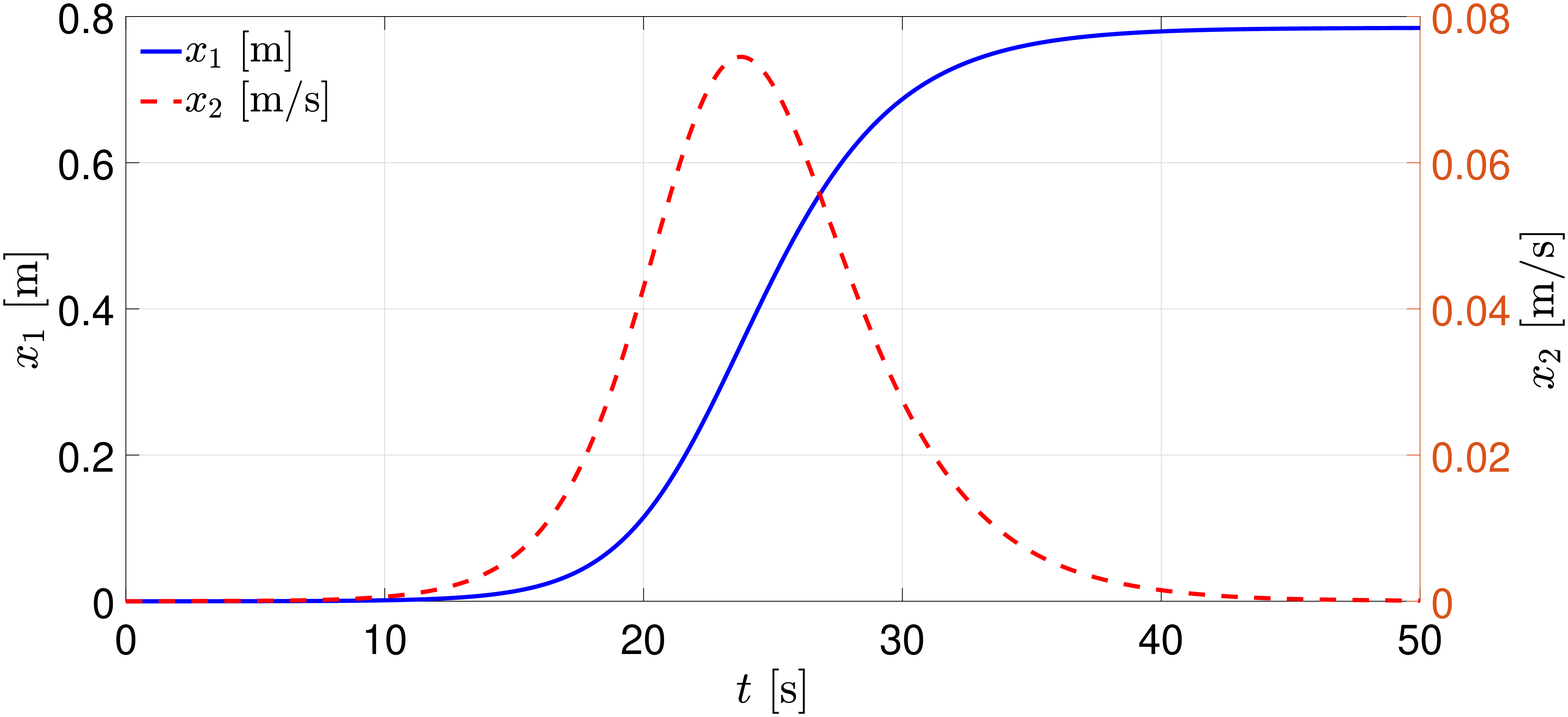}
  \caption{Earthquake phenomenon in the real-fault simulation.}
  \label{fig:earthquakesim}
\end{figure}

\begin{figure}[ht!]
  \centering 
  \includegraphics[width=8cm,height=4cm]{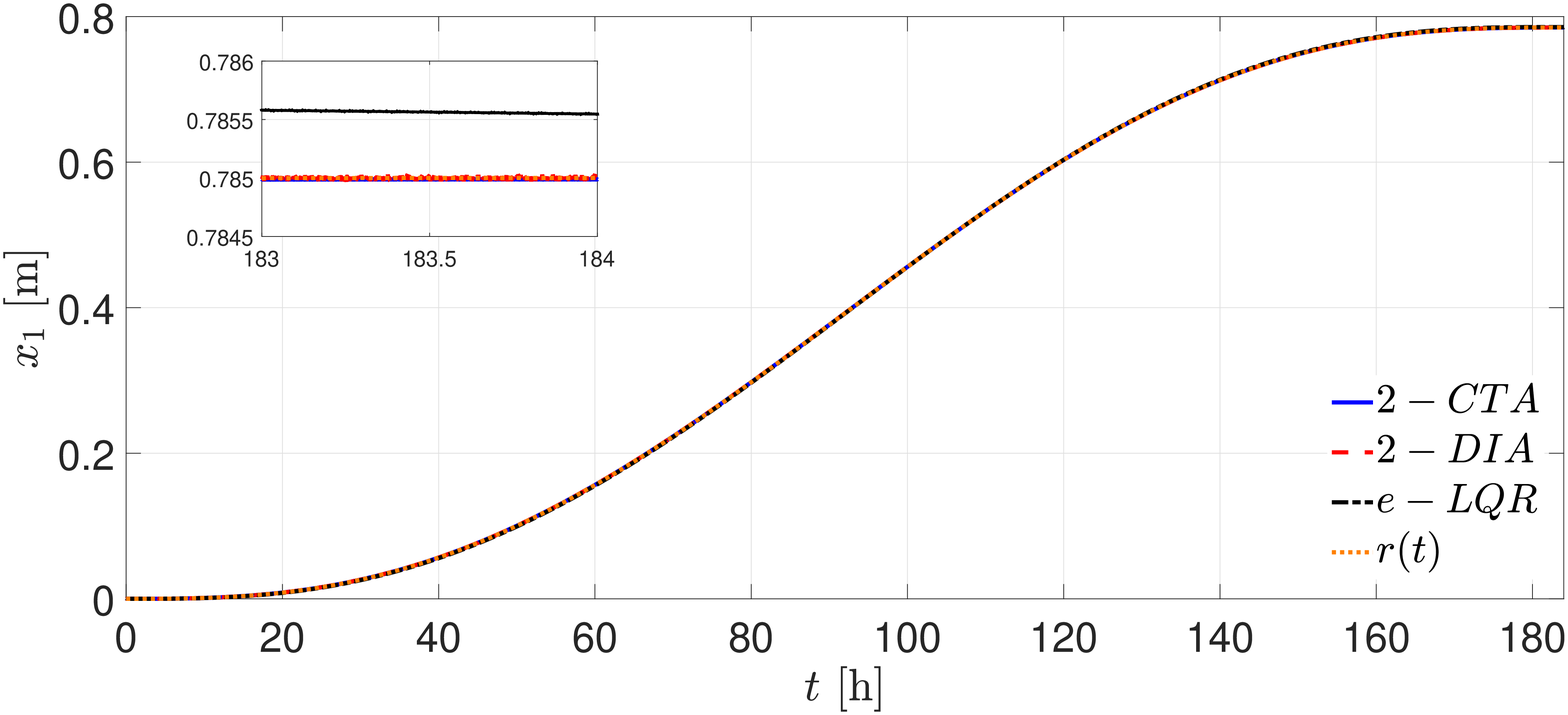}
  \includegraphics[width=8cm,height=4cm]{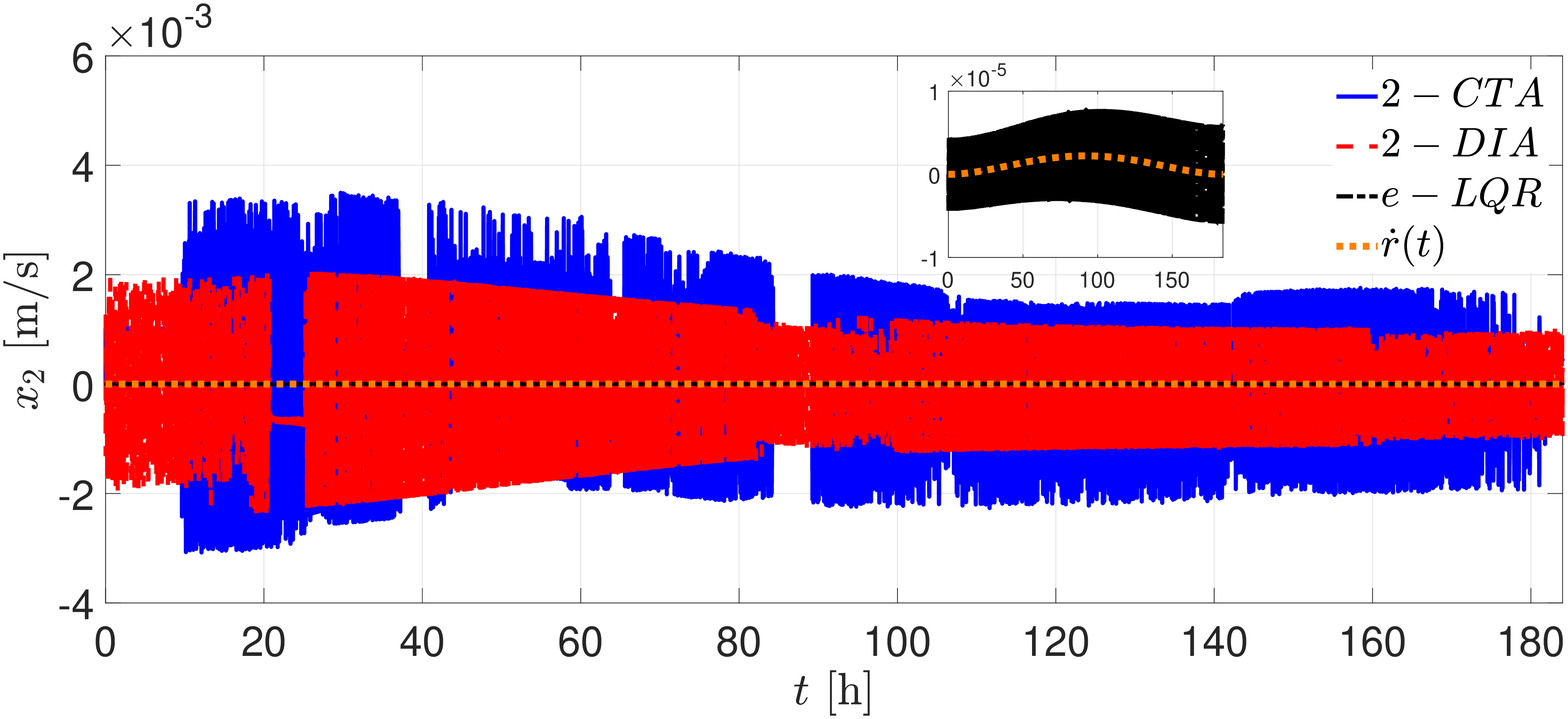}
  \caption{States in the real-fault simulations.}
  \label{fig:x1}
\end{figure}

\begin{figure}[ht!]
  \centering 
  \includegraphics[width=8cm,height=4cm]{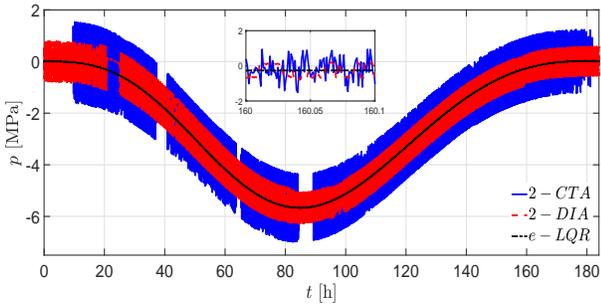}
  \caption{Control signal in the real-fault simulations.}
  \label{fig:p1}
\end{figure}

The numerical results are presented in Figs. \ref{fig:earthquakesim}-\ref{fig:p1}. In the three simulations, the slip $x_1$ is follows the desired reference in order to dissipate slowly all the stored energy avoiding earthquake-like events. This can be seen by comparing Figs. \ref{fig:earthquakesim} and \ref{fig:x1}, where the lack of a control input in the earthquake phenomenon results in the state $x_1$ evolving much faster than the controlled scenarios. Sliding-mode controllers present better results in the displacement tracking, but in terms of velocity, the e-LQR controller shows a more stable response without oscillations. Nevertheless, the three algorithms fulfil the tracking task despite the presence of the external perturbation and the use of nominal system parameters for the control design (robustness). In particular, the control signal generated by the e-LQR controller shows a smooth behaviour without demanding excess actuator response. On the other hand, the sliding-mode controllers (2-CTA and 2-DIA) show an oscillatory behaviour but not reaching high oscillations as the chattering effect (see zoom made in Fig. \ref{fig:p1}). This is because such sliding-mode algorithms make use of the discontinuous \textit{sign} function, which provokes oscillations of high (in theory infinite) frequency. Nevertheless, the final control signal that the plant experience is continuous, thanks to the type of controllers employed. The above comment could explain the presence the presence of more oscillations compared to the linear control, mainly in the velocity and in the control input.

\begin{figure*}[ht!]
  \centering 
  \includegraphics[width=5.5cm,height=4cm]{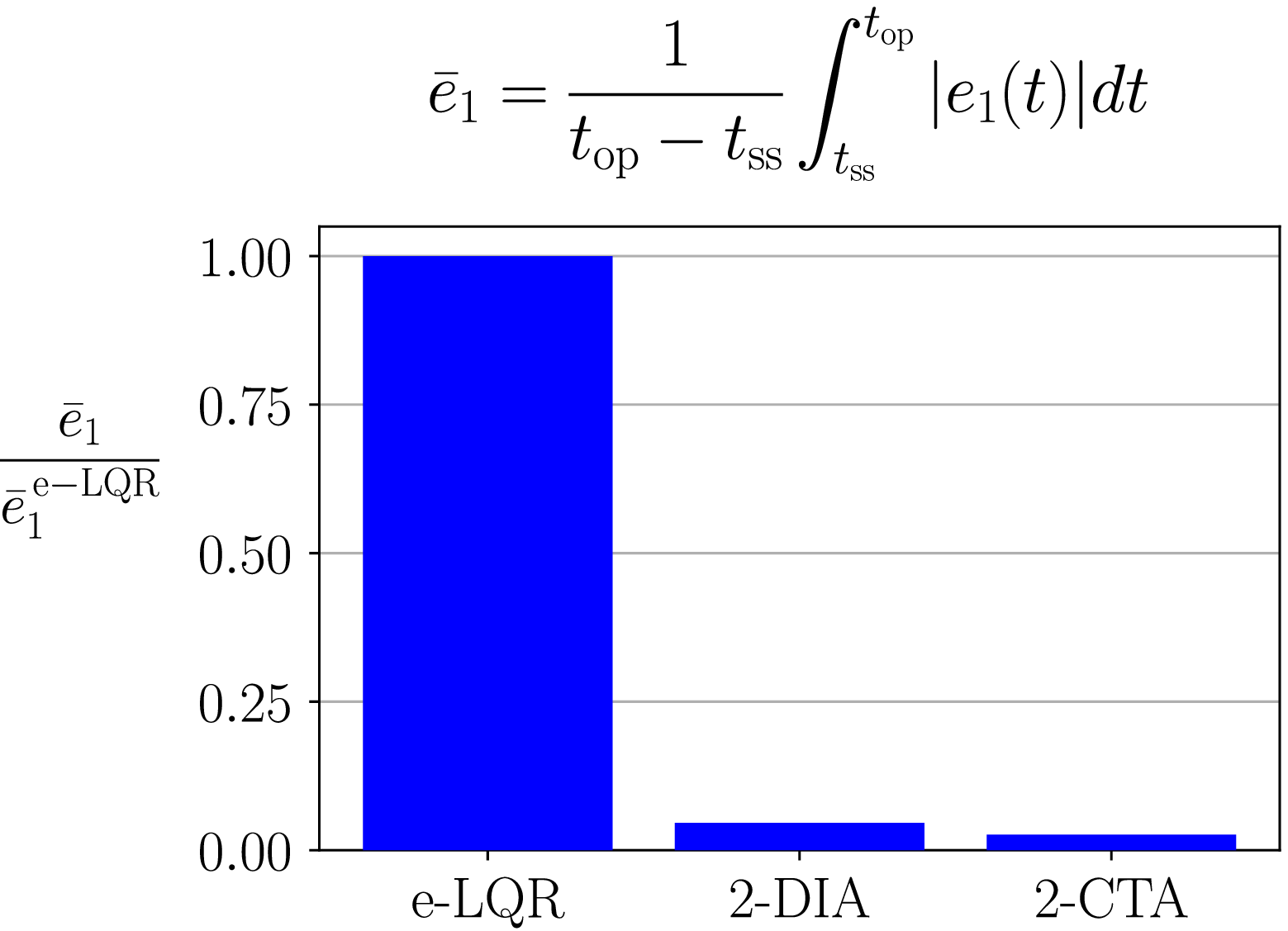}
  \includegraphics[width=5.5cm,height=3.6cm]{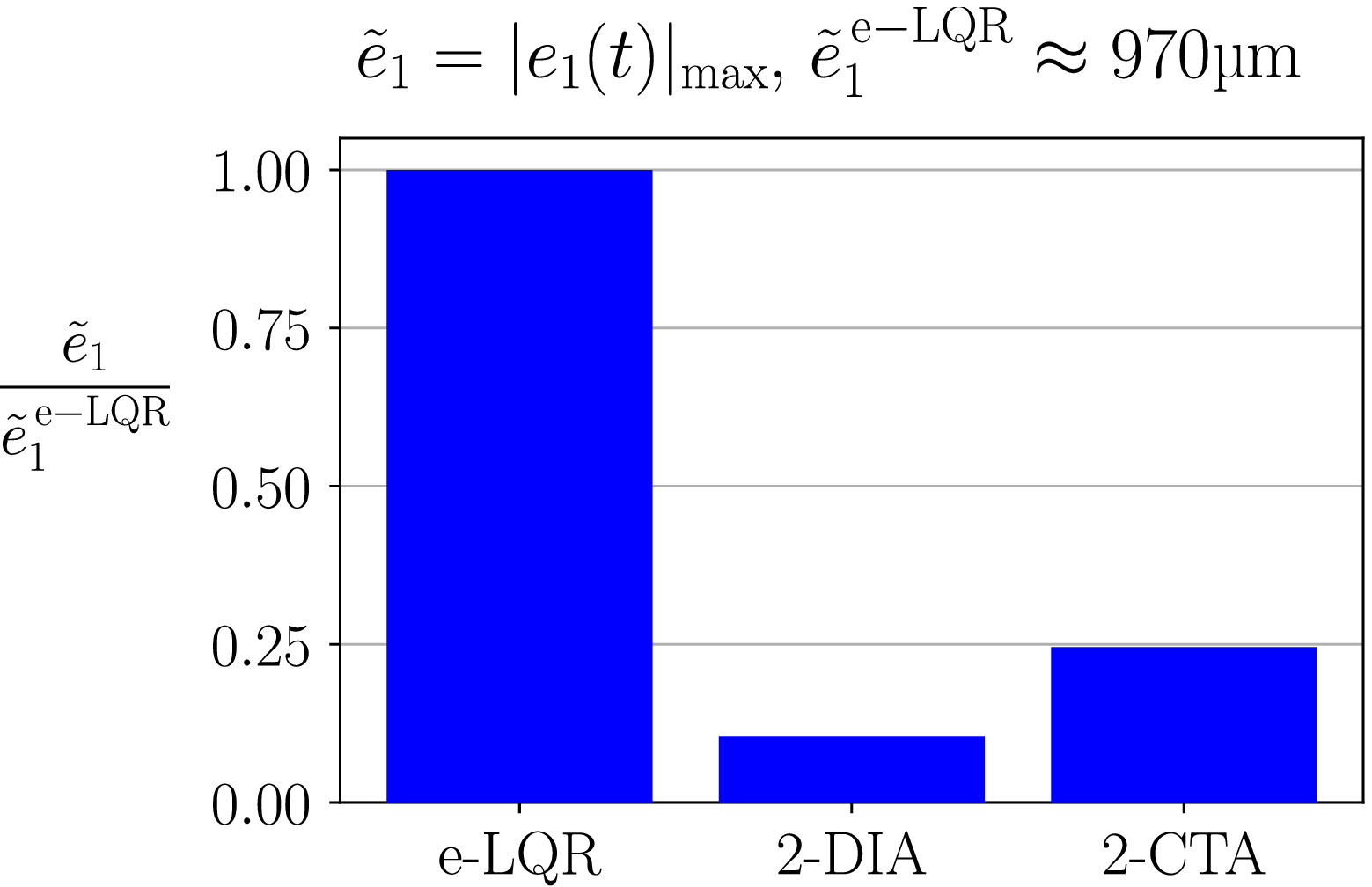}
  \includegraphics[width=5.5cm,height=4cm]{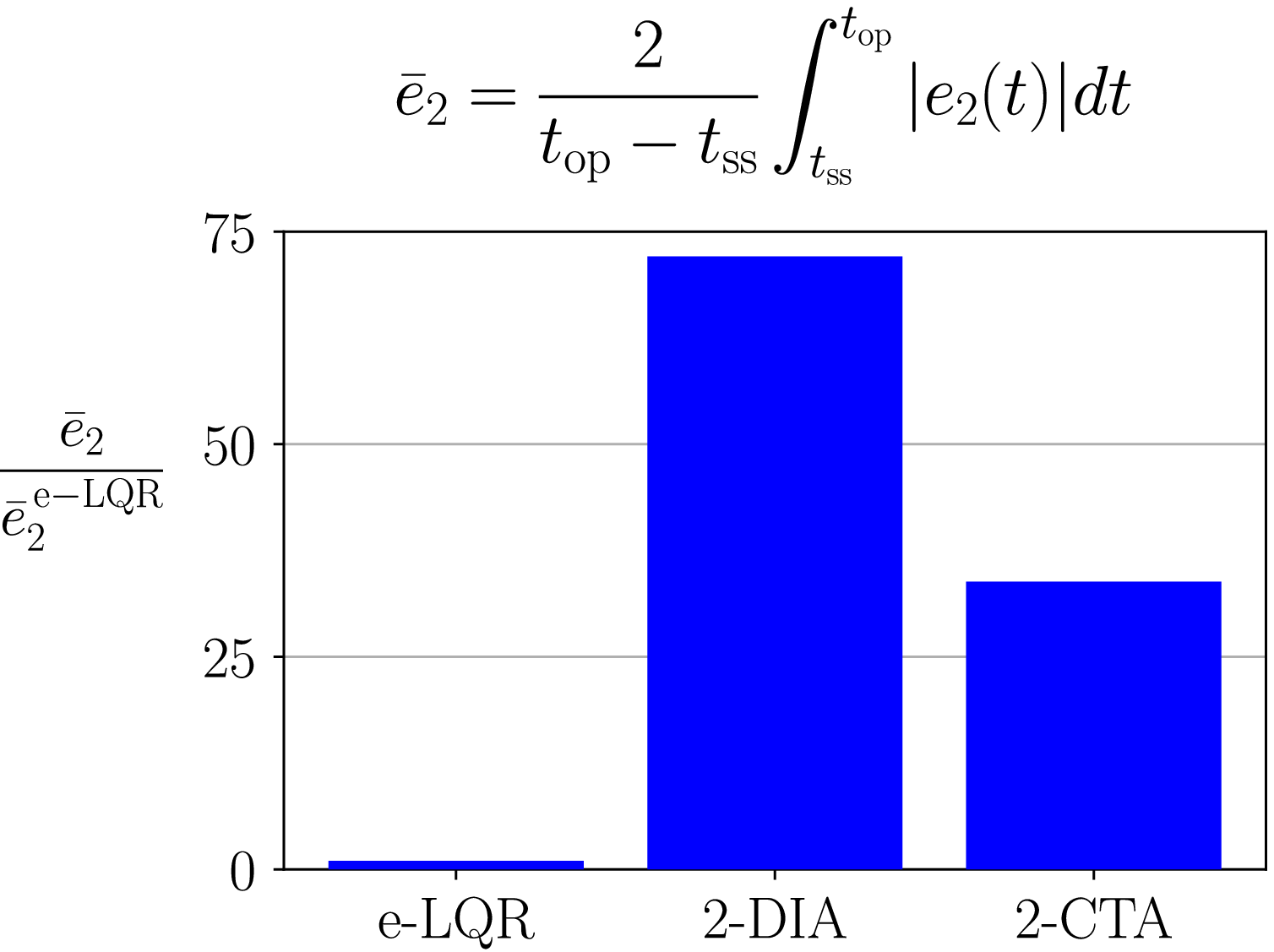}
  \includegraphics[width=5.5cm,height=3.6cm]{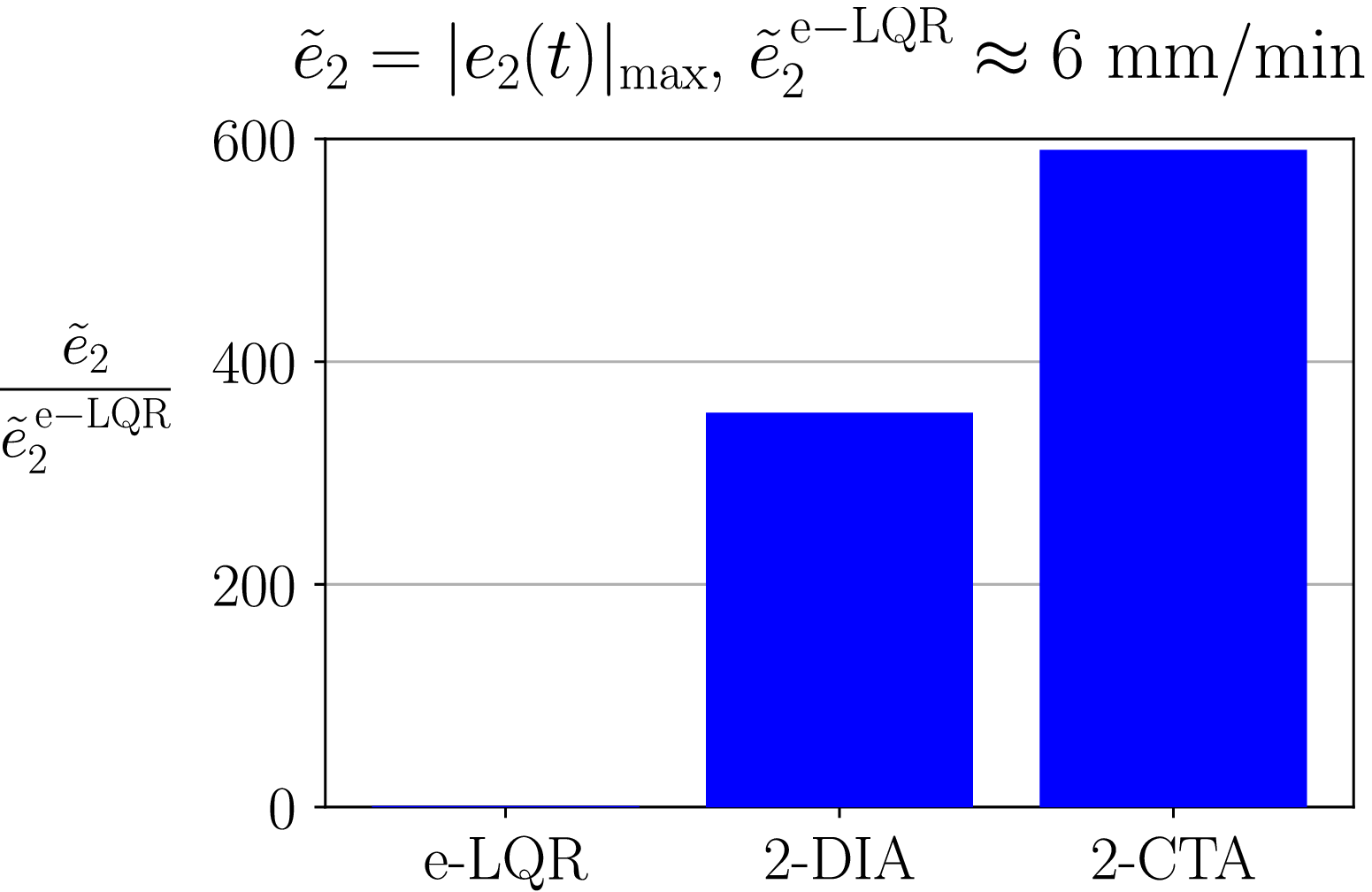}
  \includegraphics[width=5.5cm,height=4cm]{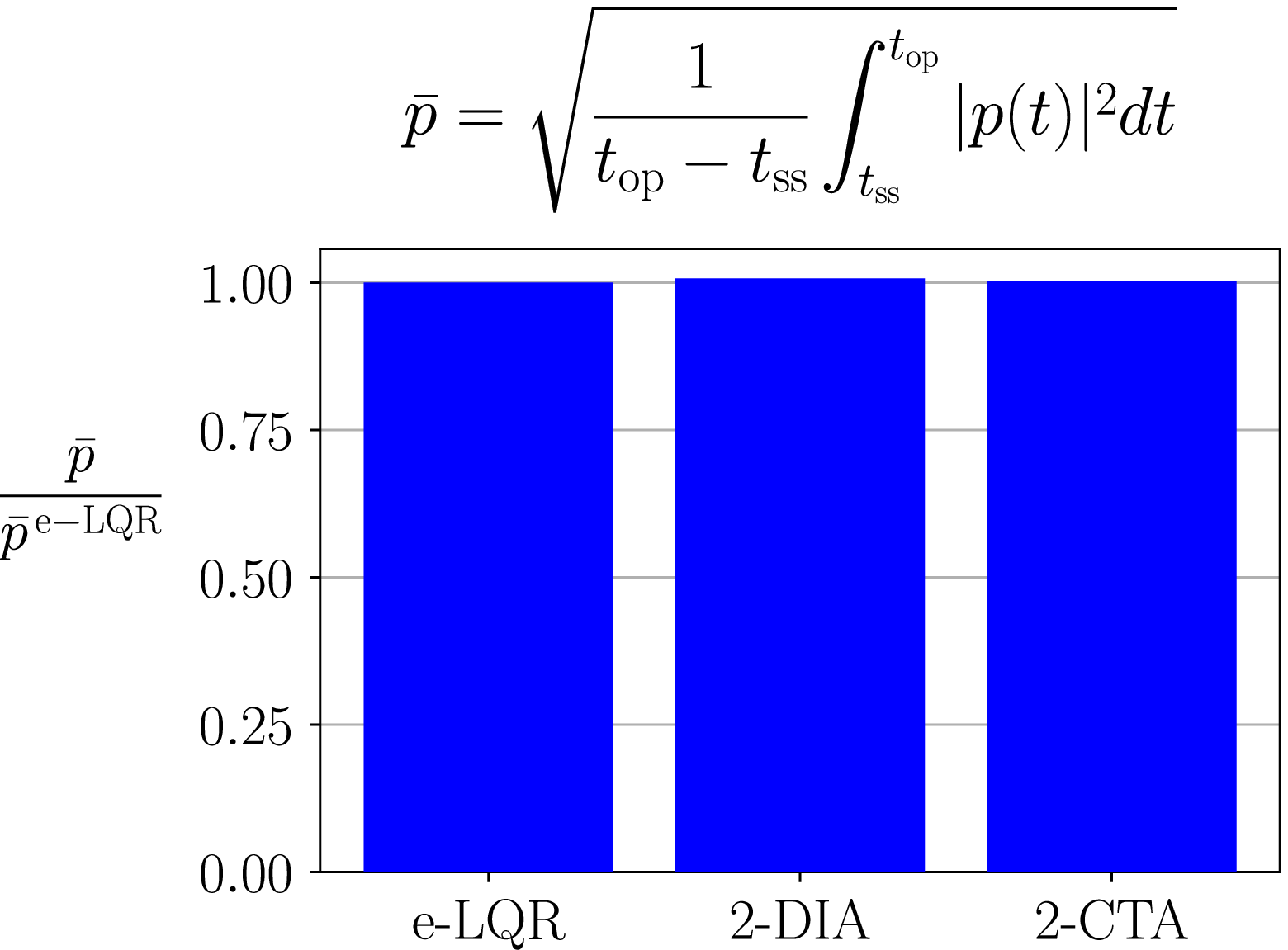}
  \caption{Errors and average power comparison in the simulations.}
  \label{fig:simindex}
\end{figure*}

To make a further comparison between the three controllers, the mean integrated error (MISE), the average power (RMS), and the maximum errors have been calculated on the steady state (a time $t=t_{ss}$). All of these were normalized w.r.t. the e-LQR results. This can be seen in Fig. \ref{fig:simindex}. The position results are smaller for the 2-DIA and the 2-CTA, but the e-LQR velocity errors are smaller. This is consistent to the results shown in latter Figs. In addition, one can observe that the maximum position error for the e-LQR control is of the order of $\sim 1$ [mm], which can be considered negligible in a total slip of $\sim 800$ [mm]. On the other hand, the maximum velocity error obtained with the 2-CTA controller was of $3.5 \times 10^{-3}$ [m/s], which twenty times smaller than the maximum slip-rate of the original earthquake event (see Fig. \ref{fig:earthquakesim}). Furthermore, the three algorithms present the same RMS value, so they spent the same amount of energy (in average) to solve the tracking problem. 

\subsection{Experimental Confirmation}

In real applications, measurements of slip and slip-rate in system \eqref{eq:springmodel} can be obtained through microseismicity measurements and other geophysical methods (\textit{e.g.}, GPS, LIDAR and topographical measurements). For our case, the lab-fault system depicted in Fig. \ref{fig:spring-slider-lab} measures the slip, ${x}_{1}$, using two LVDTs. In order to obtain the slip rate velocity $x_2$ required by the designed controllers, a robust exact filtering differentiator \cite{b:Levant-Livne-2020} has been implemented and reads as
\begin{equation}
\begin{split}
  \dot{w}_1 &= -5 \lambda_d^{\frac{1}{5}}\Sabs{w_1}^{\frac{4}{5}}+w_2, \\
  \dot{w}_2 &= -10.03 \lambda_d^{\frac{2}{5}}\Sabs{w_1}^{\frac{3}{5}}+\hat{x}_1-x_1, \\
  \dot{\hat{x}}_1 &= -9.3 \lambda_d^{\frac{3}{5}}\Sabs{w_1}^{\frac{2}{5}}+\hat{x}_2, \\
  \dot{\hat{x}}_2 &= -4.57 \lambda_d^{\frac{4}{5}}\Sabs{w_1}^{\frac{1}{5}}+\hat{x}_3, \\
  \dot{\hat{x}}_3 &= -1.1 \lambda_d\Sabs{w_1}^{0},
\end{split}
\label{eq:dif}
\end{equation} 
with $\lambda_d=1\times 10^{-5}$. Therefore, all the designed control algorithms, $p$, have been implemented as functions of the estimated slip and slip rate, \textit{i.e.}, $p(\hat{x}_1,\hat{x}_2)$ instead of $p(x_1,x_2)$. This would reduce the possible noise from the estimations due to the LVDTs measurements.

\begin{remark}
Differentiator \eqref{eq:dif} provides the second derivative of the input $x_1(t)$ while filtering the signal with a second order filter, if $\abs{{x}_1^{(3)}(t)} \leq L_d$ and $\lambda_d > L_d$ (see proof in \cite{b:Levant-Livne-2020}). Furthermore, the authors of the present paper are fully aware of the separation principle problem with nonlinear systems. However, the purpose of the following results consists of experimentally showing the convergence of the complete closed loop-system (Plant, control and differentiator).
\end{remark}

\begin{figure}[ht!]
  \centering 
  \includegraphics[width=8.5cm,height=4.2cm]{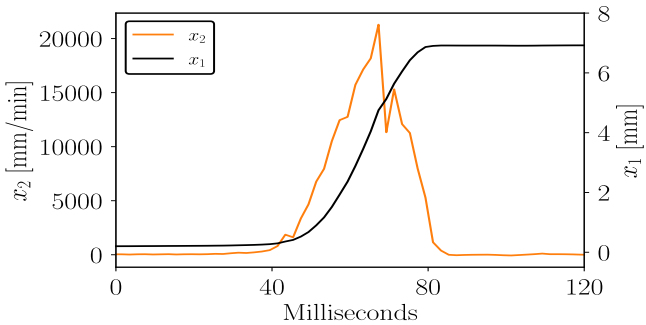}
  \caption{Earthquake phenomenon in the lab-fault experiment.}
  \label{fig:earthquakelab}
\end{figure}

The designed controllers have been tested in the laboratory apparatus shown in Fig. \ref{fig:spring-slider-lab}. The experimental process is as follows. First, the spring is compressed by a press, taking the system close to its unstable equilibrium point. Without control, the stored elastic energy of the spring is released abruptly, causing an earthquake-like response (fast seismic slip), as shown in Fig. \ref{fig:earthquakelab}.

Then, the controller \eqref{eq:ps} with $\mu_0=\mu_{res}^L, \hat{N}_0=A^L/m^L$, $\nu$ chosen as \eqref{eq:cta} for the 2-CTA and \eqref{eq:dia} for the 2-DIA, and the e-LQR control defined by \eqref{eq:u} have been implemented\footnote{Superscript $L$ means lab-fault parameters from Table \ref{tab:para}.}. The lab-fault system parameters shown in Table \ref{tab:para} are considered as nominal values of the plant for calculating the controller gains as
\begin{itemize}
  \item 2-CTA: $k_1=7.5$, $k_2=5$, $k_3=1.66 \times 10^{-4}$, $k_4=7.92 \times 10^{-5}$ and $\lambda=500$.
  \item 2-DIA: $k_{I1}=2$, $k_{I2}=5$, $k_{I3}=1.44 \times 10^{-4}$, $k_{I4}=0$ and $\lambda=500$.
  \item e-LQR: $k_1=4.17\times 10^{8}$, $k_2=6.94\times 10^{5}$, $k_3=4.17\times 10^{7}$ and $k_4=1.39\times 10^{5}$.
  \label{expgains}
\end{itemize}

\begin{figure}[ht!]
  \centering 
  \hspace{15pt}\includegraphics[width=7.3cm,height=4cm]{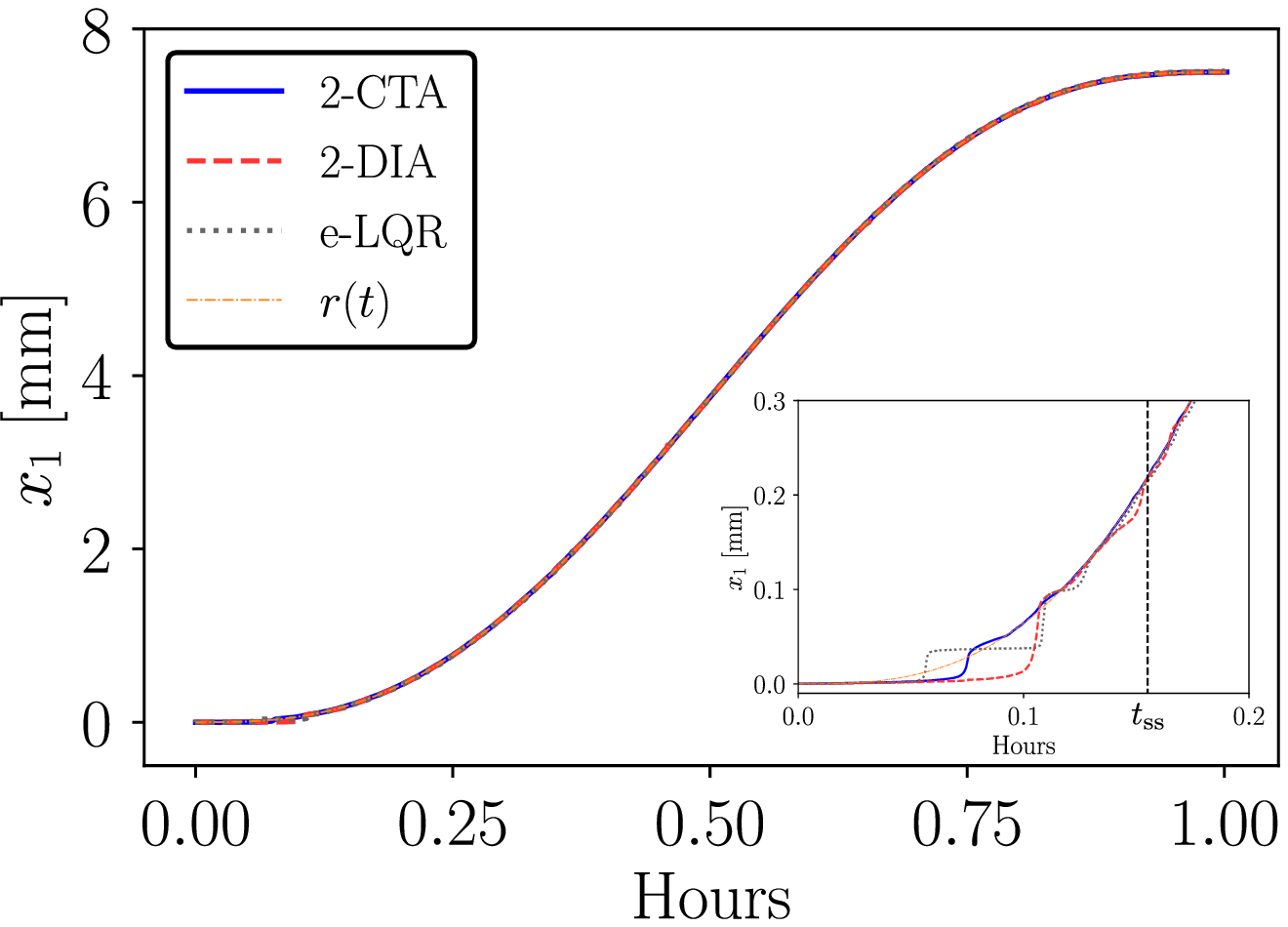}
  \includegraphics[width=7.9cm,height=4cm]{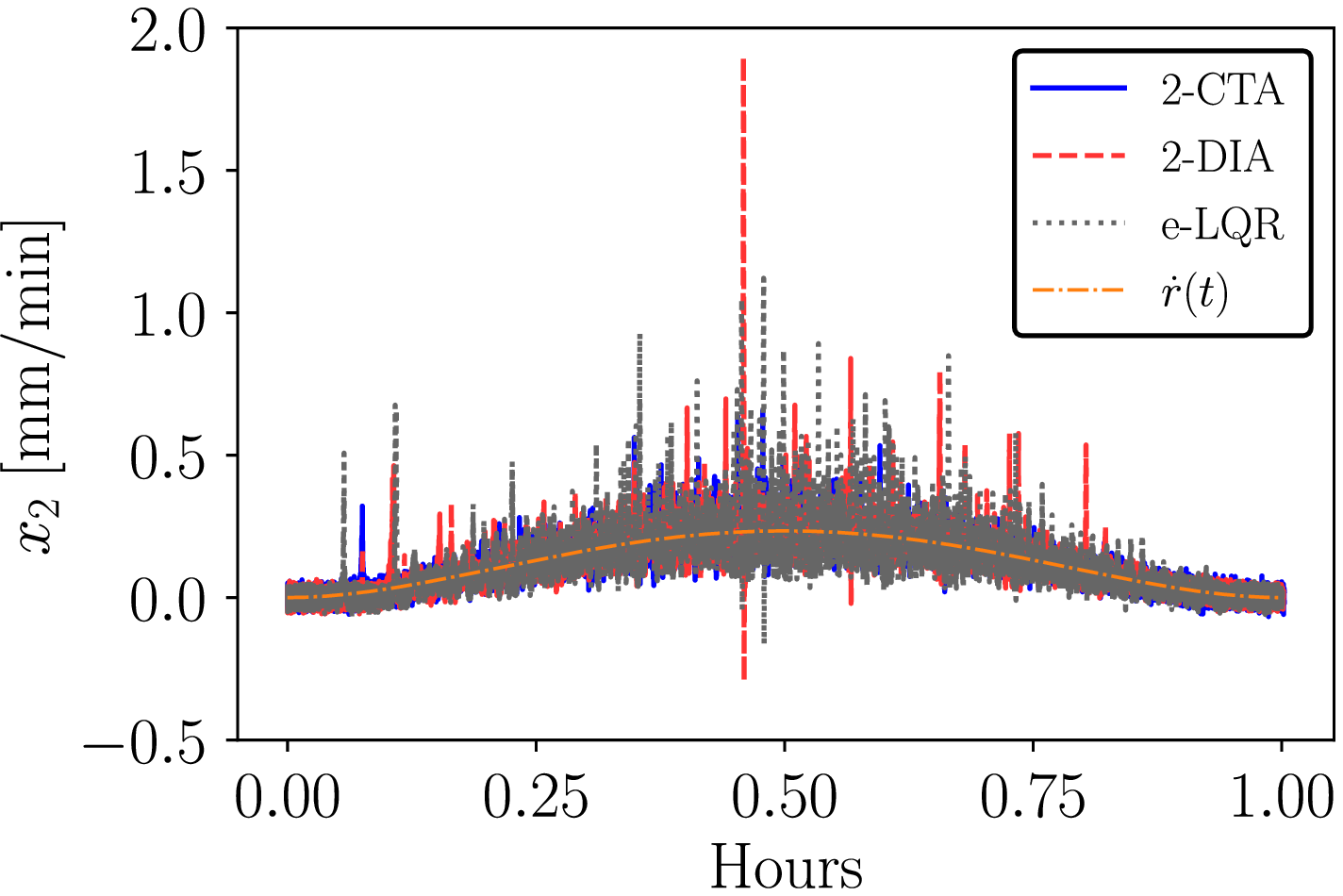}
  \caption{States in the experimental confirmation.}
  \label{fig:xe}
\end{figure}

\begin{figure}[ht!]
  \centering 
  \includegraphics[width=7.8cm,height=4.1cm]{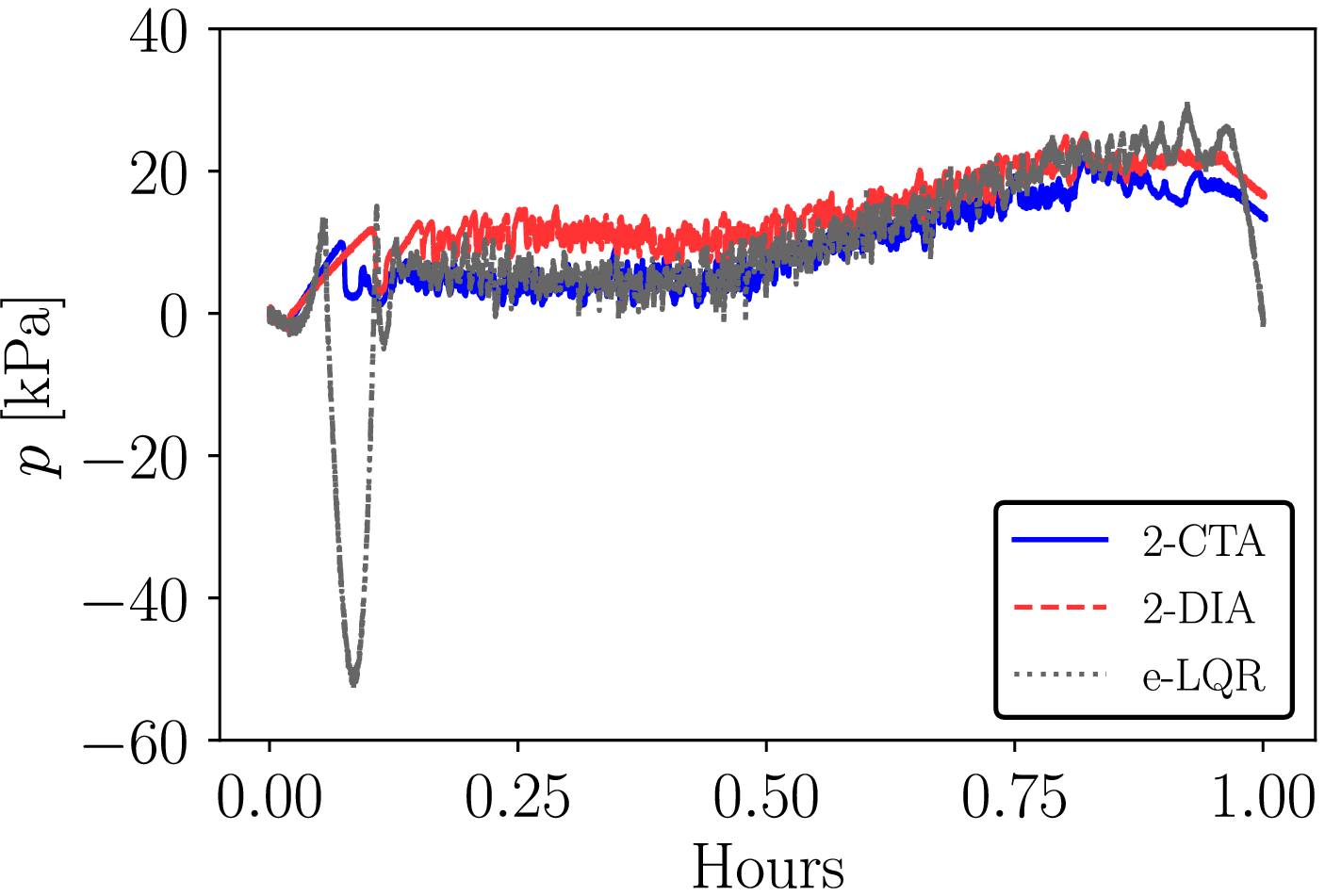} \\ \vspace{3pt}
  \hspace{6pt}\includegraphics[width=7.7cm,height=4.1cm]{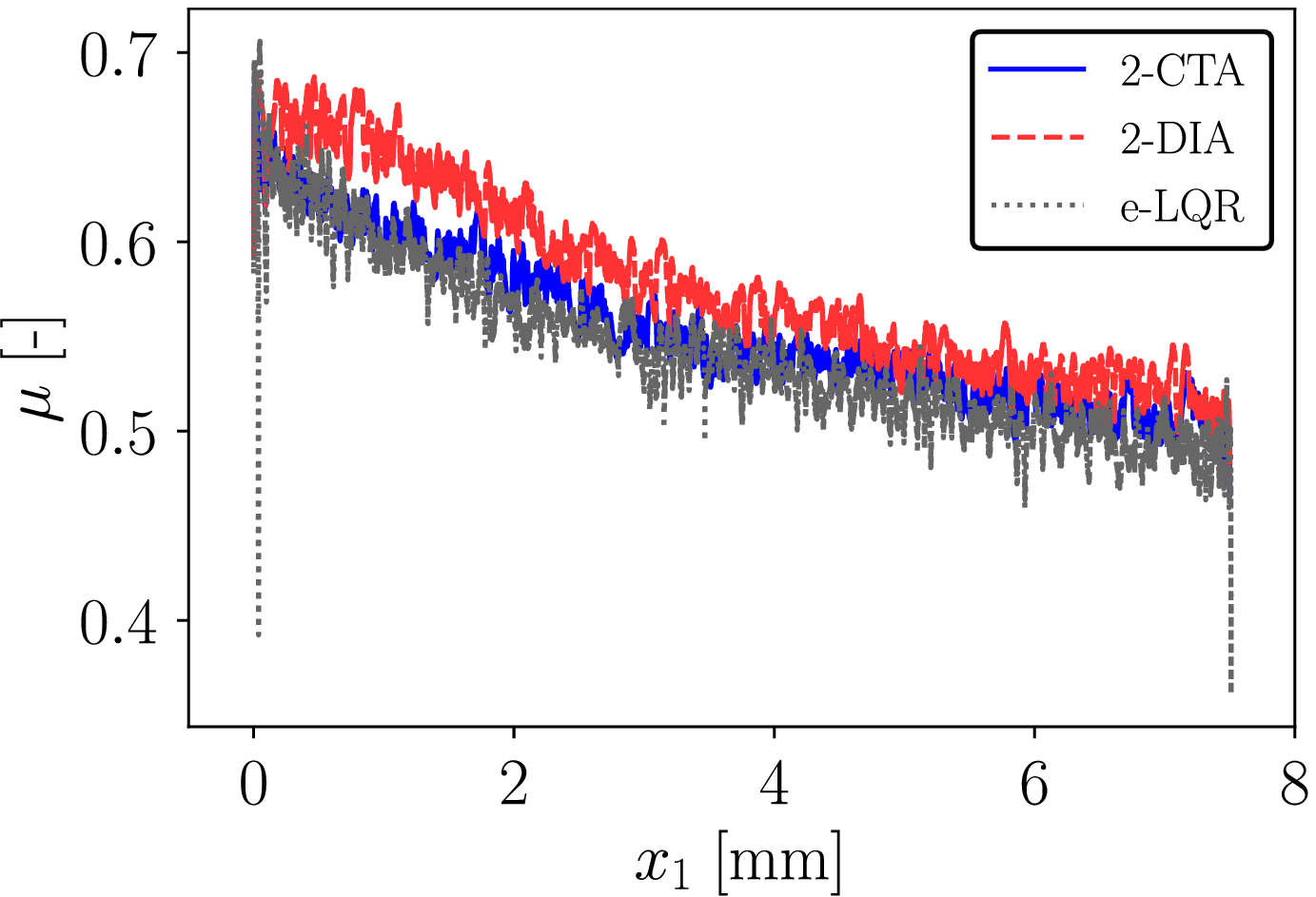}
  \caption{Control signal and friction coefficient in the experimental confirmation.}
  \label{fig:pe}
\end{figure}

The reference signal to be tracked is chosen as \eqref{eq:ref}, with $d_{max}=7.5$ [mm] and $t_{op}=1$ [h]. In other words, the controllers not only must avoid the unstable behaviour but to follow an aseismic slip equal to $7.5$ [mm] in one hour, evolving according to the sigmoid function \eqref{eq:ref}. The results are displayed in Figs. \ref{fig:xe}-\ref{fig:pe}. The three controllers are able to fulfil the task of controlling an earthquake-like instability by increasing the response time of the system (from around $80$ [ms] to $1$ [h]). The three controllers achieve the steady state around 9 [min] and they keep the states around the reference signal. The e-LQR presents a large overshoot in the control signal at the beginning of the experiment (see Fig. \ref{fig:pe}) and the 2-DIA presents higher peaks on the velocity. Fig. \ref{fig:pe} shows the three continuous control signals used for the tracking and how the measured friction coefficient evolves in time. Note how this latter is always bounded and always higher than the minimum value used for designing the controller ($\mu_{res}=0.4$). It also shows high frequency oscillations compared with the nominal friction of Fig. \eqref{fig:mu}, used in the numerical simulations. This behaviour is due to the grain size of the interfaces and other unmodelled dynamics, that the controllers successfully compensate. 

\begin{figure*}[ht!]
  \centering 
  \includegraphics[width=5.5cm,height=4cm]{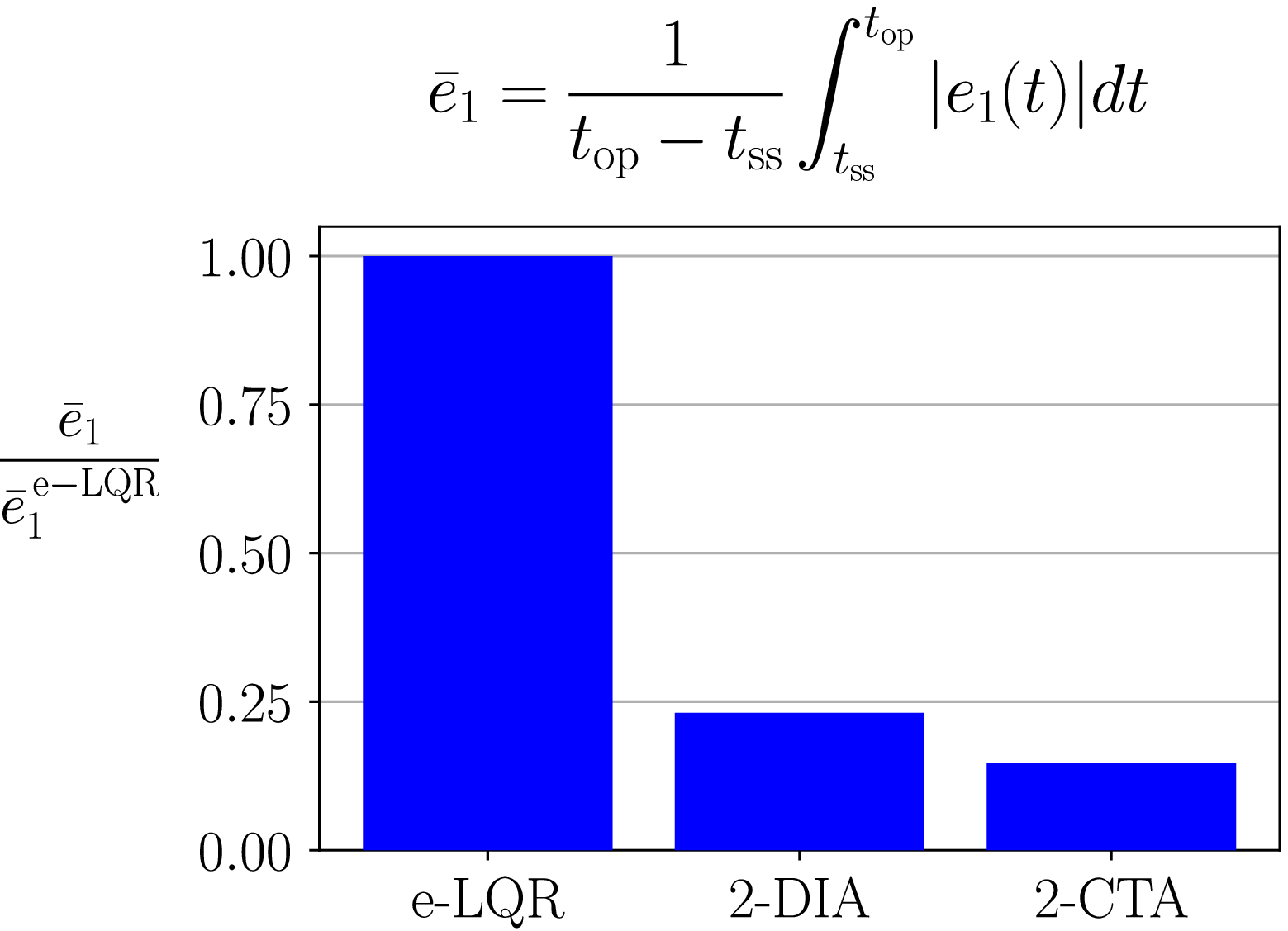}
  \includegraphics[width=5.5cm,height=3.6cm]{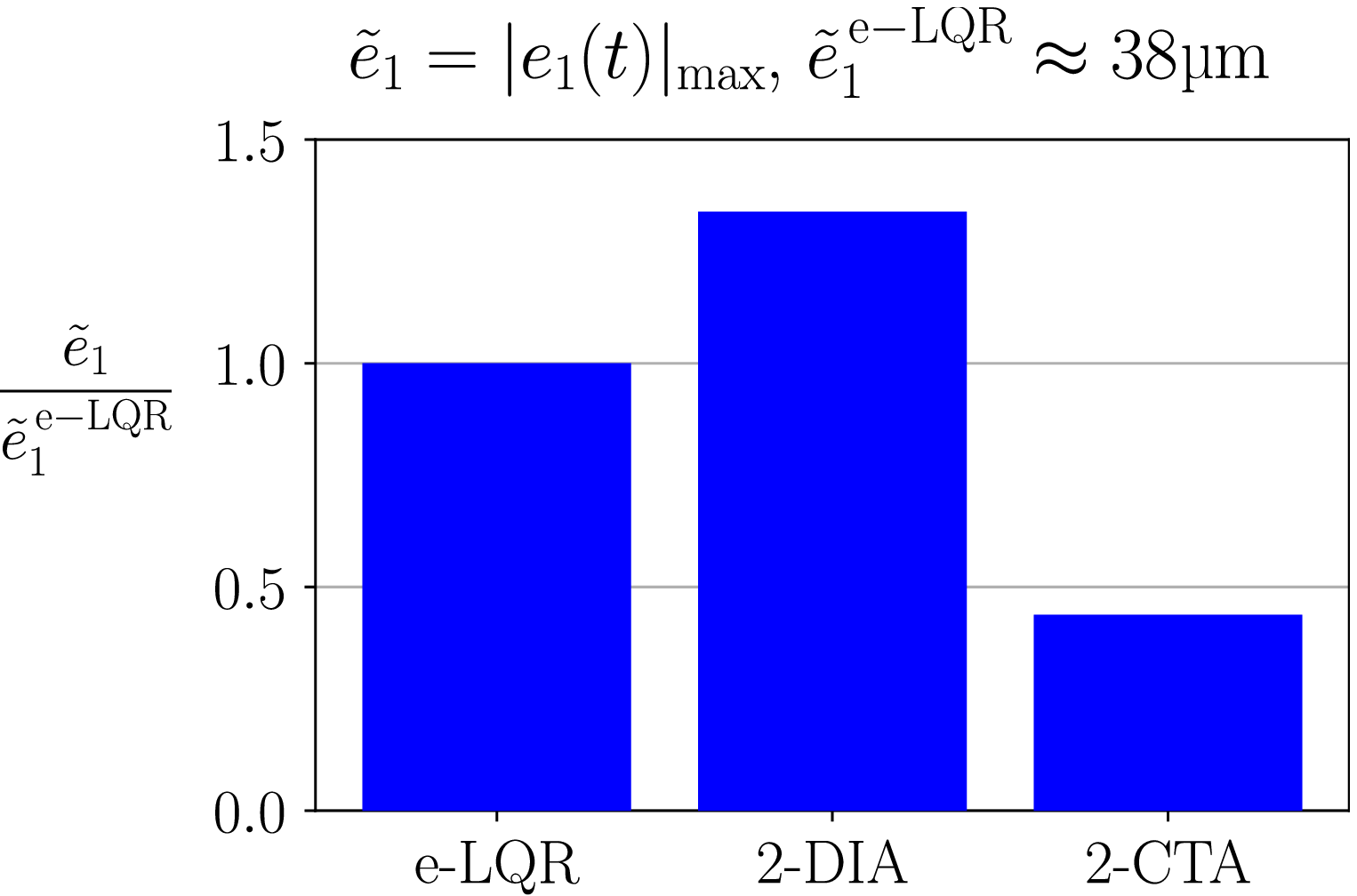}
  \includegraphics[width=5.5cm,height=4cm]{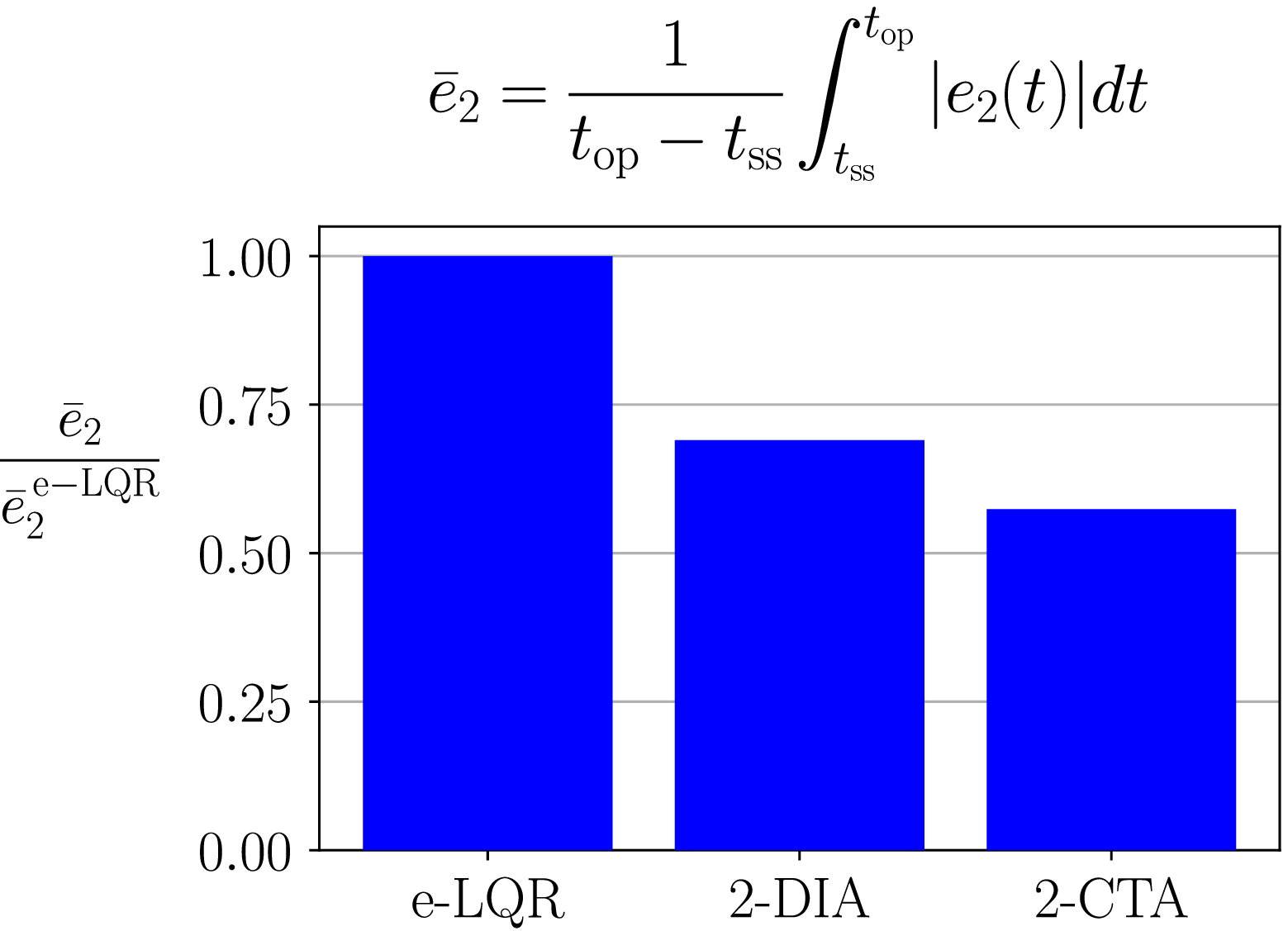}
  \includegraphics[width=5.5cm,height=3.6cm]{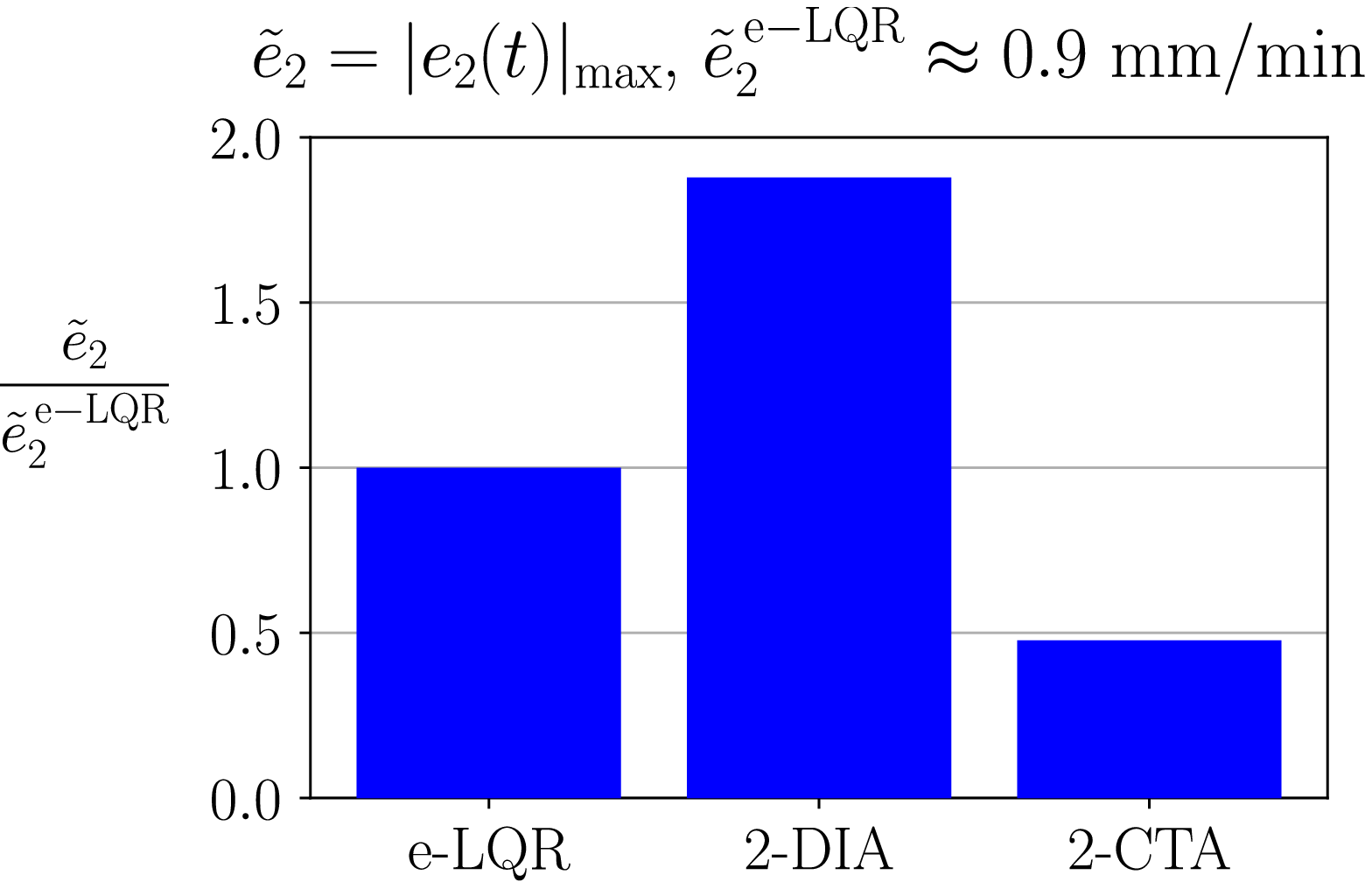}
  \includegraphics[width=5.5cm,height=4cm]{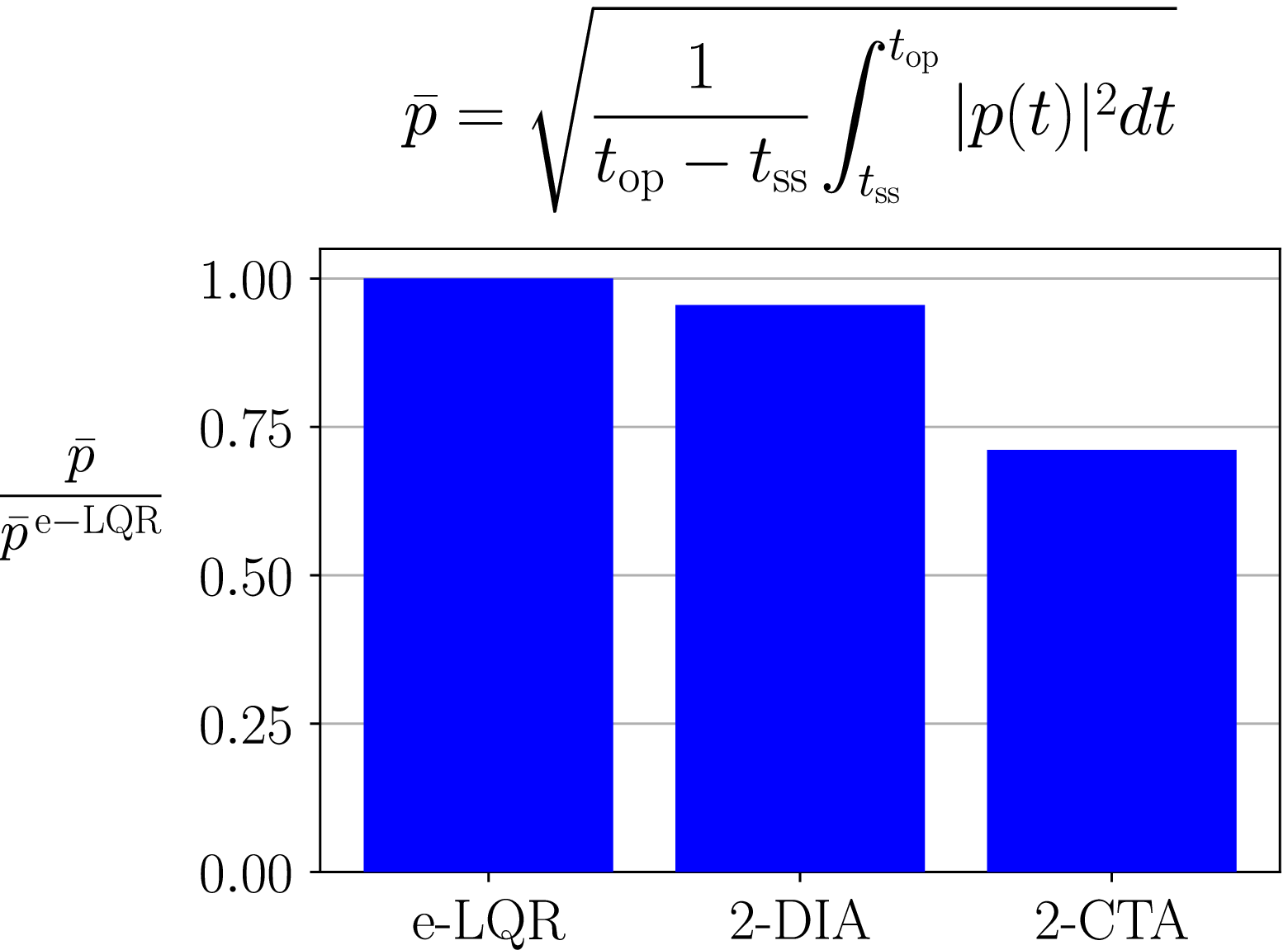}
  \caption{Errors and average power comparison in the experiments.}
  \label{fig:labindex}
\end{figure*}

Assuring a fair comparison between the three presented controllers is not straightforward, due to the presence of different uncertainties and disturbances in each experiment. In particular, the tested samples are not the same between different tests, the experiments do not initiate exactly in the same initial point, and better gains may exist to optimally tune the three different controllers. Furthermore, as discussed before, both control strategies present different theoretical properties. Nevertheless, the error comparison used in the simulations has been made again for the experimental results. This is presented in Fig. \ref{fig:labindex}. The lowest errors were obtained with the 2-CTA using the smallest average power, but the three algorithms present negligible absolute errors. A video of the experiment is available in \cite{b:CoQuake_video}. 

\section{Conclusions}
\label{sec:Conclusions}

In this paper, a challenging, emerging application of robust nonlinear control theory is presented for preventing natural and anthropogenic seismicity. The control is designed to perform tracking of a slow reference based on a reduced-order model for earthquakes. Two types of controllers are presented: the first one is based on sliding-mode theory and the other on LQR control. The first one results in local finite-time convergence of the tracking error, while the second presents global exponential stability. Both controllers are designed to generate a continuous control signal and use integral action to compensate different kinds of perturbations. The algorithms have been tested and compared in numerical simulations over a real-fault system and in a specially designed experimental apparatus, showing that both types of controllers succeed in achieving tracking to a new stable equilibrium of lower energy. In the simulations, the best precision for the slip error was obtained with the sliding-mode algorithms, but the LQR control was better achieving a smaller error in the slip-rate error. With respect to the experimental tests, both controllers were successful in compensating unmodelled dynamics and parameter uncertainties present in the real systems, but the 2-CTA sliding-mode algorithm achieved the best results. As far as it concerns, the performance of the controllers could be upscaled to a real earthquake due to the existence of scaling laws between both faults. The design of controllers based on more detailed and complete models representing the earthquake phenomenon (\textit{e.g.}, coupled elastic and diffusion partial differential equations) remains as future work.

\appendices

\section*{Acknowledgment}

The authors would like to acknowledge the support of the European Research Council (ERC) under the European Union’s Horizon 2020 research and innovation program (Grant agreement no. 757848 CoQuake).

\section{Gain scaling}
\label{sec:app1}

To show how the experimental results can be upscaled to a more realistic earthquake event, the lab-fault system parameters were scaled to obtained the real-fault parameters, according to Table \ref{tab:para}. Furthermore, the controller gains used in the experiments \eqref{expgains} were manipulated in the same sense to obtain the set of gains \eqref{simgains} for the numerical simulations.

Recalling \eqref{eq:ps}, \eqref{eq:cta}, \eqref{eq:dia} and \eqref{eq:u}, the controllers to be used in the lab-fault experiments are defined as
\begin{equation*}
	\begin{split}
	\textup{2-CTA:}\quad p^L = & \left[ - \lambda^{\frac{2}{3}} k_1\Sabs{e_1^L}^{\frac{1}{3}}- \lambda^{\frac{1}{2}} k_2\Sabs{e_2^L}^{\frac{1}{2}} + \xi_1^L \right] \frac{1}{\mu_0^L \hat{N}_0^L},\\
	\dot{\xi}_1^L = & - \lambda k_3\Sabs{e_1^L}^{0}- \lambda k_4\Sabs{e_2^L}^{0},
	\end{split}
\end{equation*}
\begin{equation*}
	\begin{split}
	\textup{2-DIA:}\quad p^L = & - \left[ \lambda^{\frac{1}{2}} k_{I2}\Sabs{\Sabs{e_2^L}^{\frac{3}{2}} + \lambda^{\frac{1}{2}} k_{I1}^{\frac{3}{2}} e_1^L}^{\frac{1}{3}} + \xi_1^L \right] \frac{1}{\mu_0^L \hat{N}_0^L},\\
	\dot{\xi_1}^L = & - \lambda k_{I3} \Sabs{e_1^L + \lambda^{-\frac{1}{2}} k_{I4} \Sabs{e_2^L}^{\frac{3}{2}}}^{0},
	\end{split}
\end{equation*}
\begin{flalign*}
&\textup{e-LQR:}\quad p^L = -k_1 x_1^L -k_2 x_2^L -k_3 \xi_1^L -k_4 \xi_2^L,&
\end{flalign*}
using the superscript $L$ as notation.

According to Table \ref{tab:para}, the upscaled pressure, slip and slip-rate, error variables and integral terms can be obtained as $p^R=p^L \lambda_p$, $e_1^R=e_1^L \lambda_\delta$, $x_1^R=x_1^L \lambda_\delta$, $e_2^R=e_2^L \lambda_v$, $x_2^R=x_2^L \lambda_v$, $\xi_1^R=\xi_1^L \lambda_\delta \lambda_t$, $\xi_2^R=\xi_2^L \lambda_\delta \lambda_t^2$, $\lambda_v=\nicefrac{\lambda_\delta}{\lambda_t}$. Therefore, the latter controllers for the real-fault simulations read as
\begin{equation*}
	\begin{split}
	\textup{2-CTA:}\quad  p^R = & \left[ - \lambda^{\frac{2}{3}} \bar{k}_1\Sabs{e_1^R}^{\frac{1}{3}}- \lambda^{\frac{1}{2}} \bar{k}_2\Sabs{e_2^R}^{\frac{1}{2}} + \xi_1^R \right] \frac{1}{\mu_0^L \hat{N}_0^L},\\
	\dot{\xi}_1^R = & - \lambda \bar{k}_3\Sabs{e_1^R}^{0}- \lambda \bar{k}_4\Sabs{e_2^R}^{0},
	\end{split}
\end{equation*}
where $\left(\bar{k}_{1},\,\bar{k}_{2},\,\bar{k}_{3},\,\bar{k}_{4}\right)=\left(k_{1}\frac{\lambda_p}{\lambda_\delta^{\nicefrac{1}{3}}},\,k_{2}\frac{\lambda_p}{\lambda_v^{\nicefrac{1}{2}}},\,k_{3}\frac{\lambda_p}{\lambda_t},\,k_{4}\frac{\lambda_p}{\lambda_t} \right)$,
{\small
\begin{equation*}
	\begin{split}
	\textup{2-DIA:}\quad p^R = & - \left[ \lambda^{\frac{1}{2}} \bar{k}_{I2}\Sabs{\Sabs{e_2^R}^{\frac{3}{2}} + \lambda^{\frac{1}{2}} \bar{k}_{I1}^{\frac{3}{2}} e_1^R}^{\frac{1}{3}} + \xi_1^R \right] \frac{1}{\mu_0^L \hat{N}_0^L},\\
	\dot{\xi}_1^R = & - \lambda \bar{k}_{I3} \Sabs{e_1^R + \lambda^{-\frac{1}{2}} \bar{k}_{I4} \Sabs{e_2^R}^{\frac{3}{2}}}^{0},
	\end{split}
\end{equation*}}
where $$\left(\bar{k}_{I1},\,\bar{k}_{I2},\,\bar{k}_{I3},\,\bar{k}_{I4}\right)=\left(k_{I1}\frac{\lambda_v}{\lambda_\delta^{\nicefrac{2}{3}}},\,k_{I2}\frac{\lambda_p}{\lambda_v^{\nicefrac{1}{2}}},\,k_{I3}\frac{\lambda_p}{\lambda_t},\,k_{I4}\frac{\lambda_\delta}{\lambda_v^{\nicefrac{3}{2}}} \right),$$
\begin{flalign*}
&\textup{e-LQR:}\quad p^R = -\bar{k}_1 x_1^R -\bar{k}_2 x_2^R -\bar{k}_3 \xi_1^R -\bar{k}_4 \xi_2^R,&
\end{flalign*}
where $\left(\bar{k}_{1},\,\bar{k}_{2},\,\bar{k}_{3},\,\bar{k}_{4}\right)=\left(k_{1}\frac{\lambda_p}{\lambda_\delta},\,k_{2}\frac{\lambda_p}{\lambda_v},\,k_{3}\frac{\lambda_p}{\lambda_\delta \lambda_t},\,k_{4}\frac{\lambda_p}{\lambda_\delta \lambda_t^2} \right)$, 
and using the superscript $R$ as notation.

\bibliography{Bibliografias}
\bibliographystyle{ieeetr}

\begin{IEEEbiography}[{\includegraphics[width=1in,height=1.25in,clip,keepaspectratio]{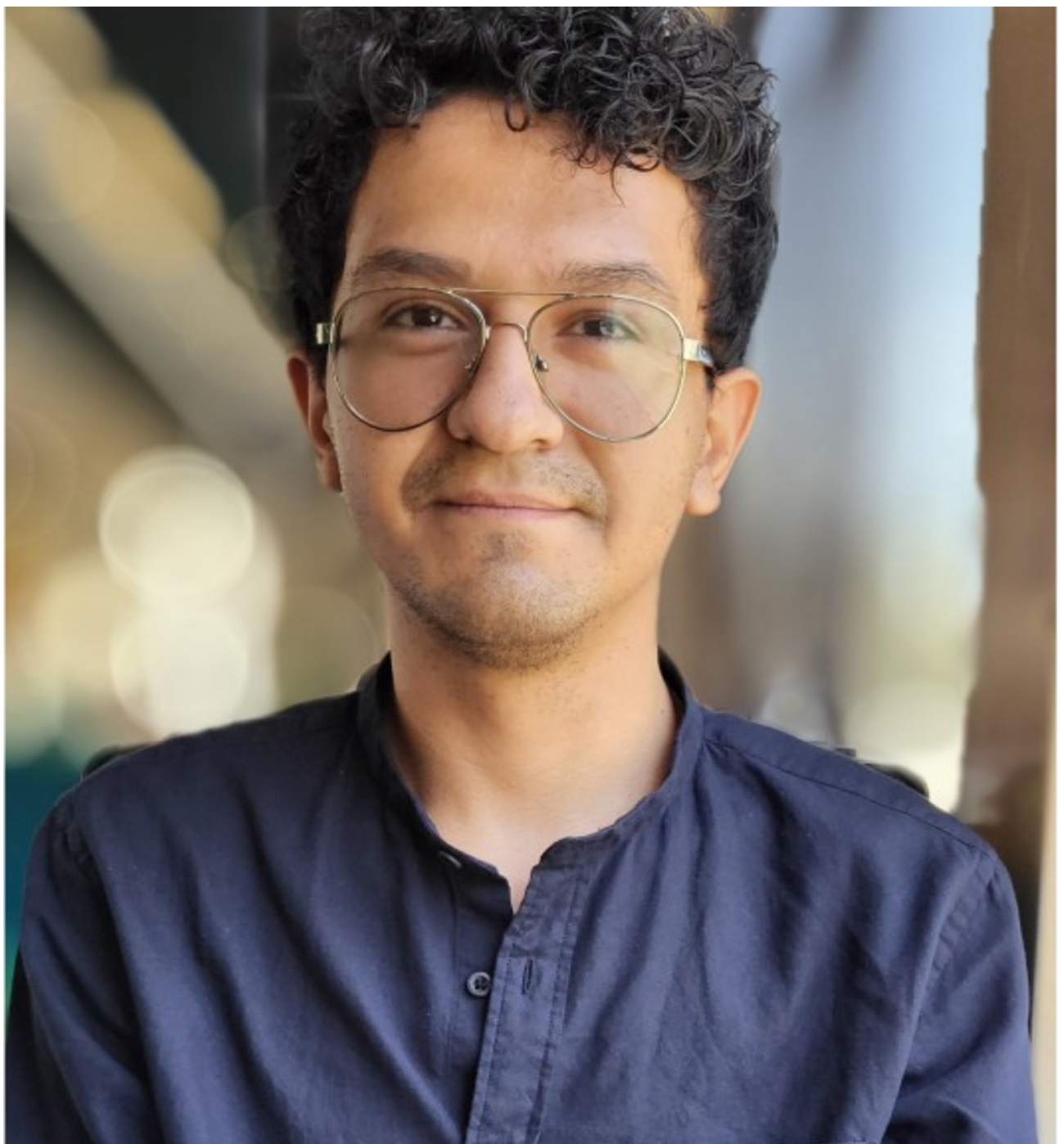}}]{Diego Guti\'errez-Oribio} was born in M\'exico and received his M.S. degree and Bachelor degree in Mechatronic Engineering from the National Autonomous University of Mexico (UNAM), Mexico City, Mexico in 2016 and 2013, respectively. In 2021, he received his Ph.D. in Electrical Engineering - Automatic Control (with honours) at UNAM. He is currently a post-doctoral member in the GeM Laboratory in the Ecole Centrale, Nantes, France. His current research of interests include the control of mechanical systems, robust and nonlinear control, continuous higher-order sliding-mode control and control of earthquake phenomena.
\end{IEEEbiography}

\begin{IEEEbiography}[{\includegraphics[width=1in,height=1.25in,clip,keepaspectratio]{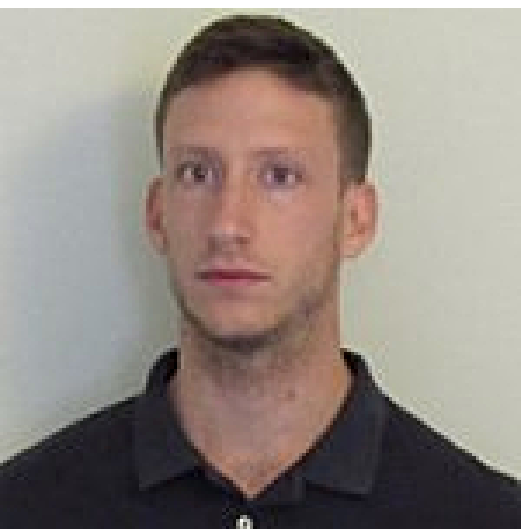}}]{Georgios Tzortzopoulos} does his PhD thesis on earthquake control using pertinent fault stimulating experimental techniques in the laboratory. He graduated in 2016 from the National Technical University of Athens (NTUA) holding a Diploma in Civil Engineering and a Master’s Degree in Structural Engineering. He did his Master Thesis at the Laboratory of Structural Analysis and Antiseismic Research of NTUA. In the meantime, he fulfilled his military obligations in Greece. His current research interests are oriented towards fault mechanics and earthquake control using advanced experimental techniques.
\end{IEEEbiography}

\begin{IEEEbiography}[{\includegraphics[width=1in,height=1.25in,clip,keepaspectratio]{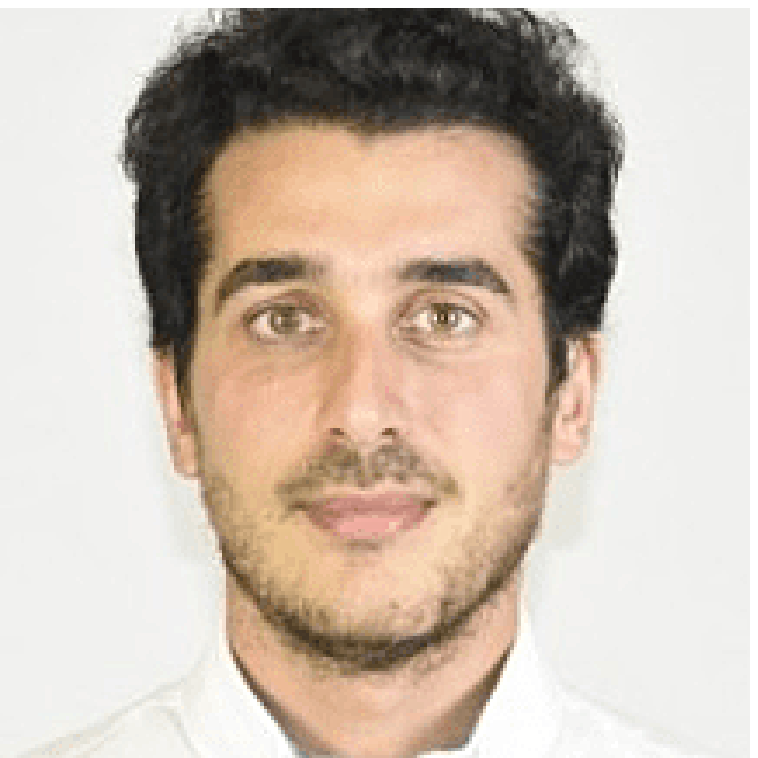}}]{Ioannis Stefanou} is Professor at the Ecole Centrale de Nantes (ECN), France. He studied civil engineering, mechanics and applied mathematics at the National Technical University of Athens. He then did his PhD thesis at the Laboratory of Geomaterials of the same institution. His main research topics are the mechanics of (geo-)materials, structural dynamics, geomechanics, fault reactivation and earthquake control. He is the PI of the ERC-StG project "Controlling earthQuakes - CoQuake" (www.coquake.eu) and of the Connect Talent project "Controlling Extreme EVents - CEEV" (www.blastructures.eu) awarded by the Pays de la Loire.
\end{IEEEbiography}

\begin{IEEEbiography}[{\includegraphics[width=1in,height=1.25in,clip,keepaspectratio]{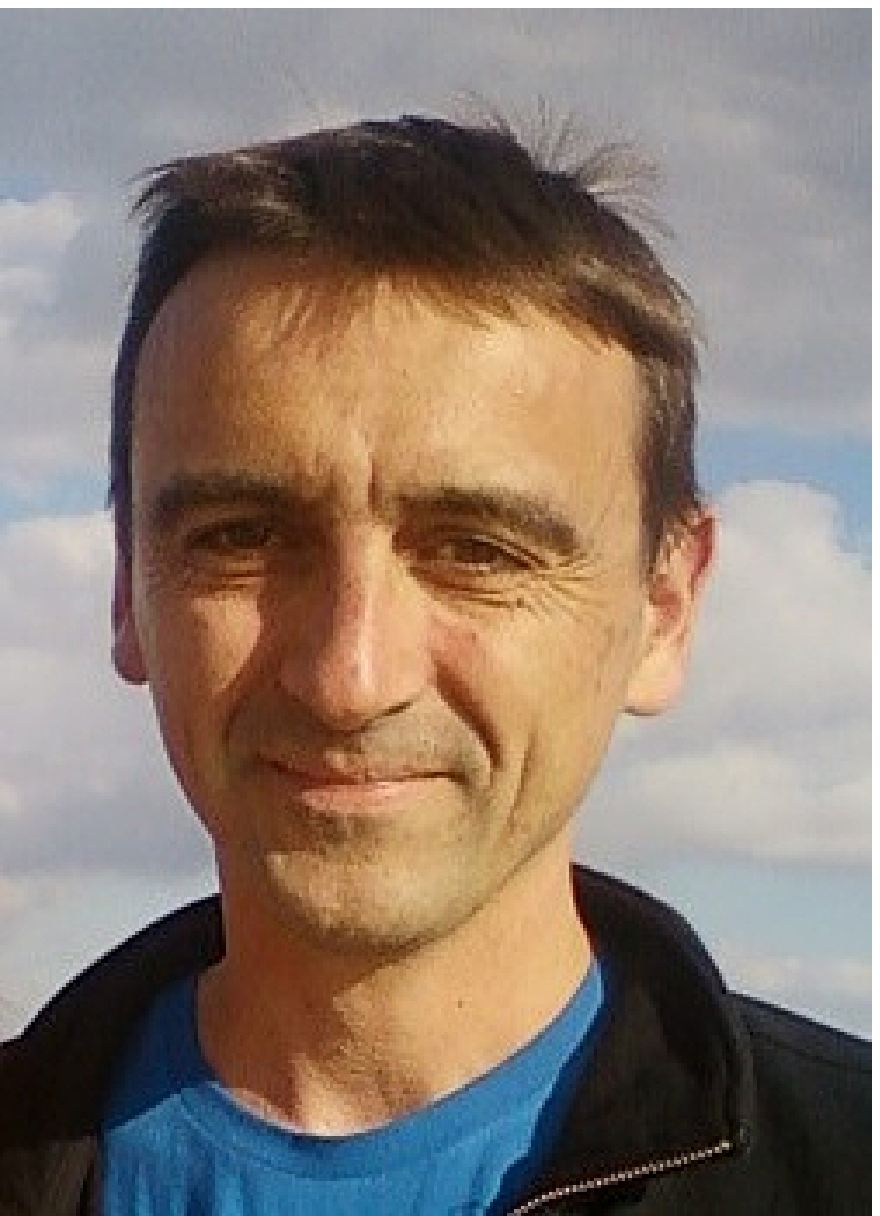}}]{Franck Plestan} received the Ph.D. in Automatic Control from the Ecole Centrale de Nantes, France, in 1995. From September 1996 to August 2000, he was with Louis Pasteur University, Strasbourg, France. In September 2000, he joined the Ecole Centrale de Nantes, Nantes, France where he is currently Professor. His research interests include robust control (adaptive/higher order sliding mode, time-delay systems) and nonlinear observer design. He is also working in several application domains as pneumatic actuators,  automotive, flying systems, renewable energy systems. 
\end{IEEEbiography}

\end{document}